\documentclass[11pt,reqno]{amsart}
\usepackage{amssymb}
\usepackage{amscd}
\usepackage{amsxtra}
\usepackage[mathscr]{eucal}

\theoremstyle{plain}
\newtheorem{theorem}{Theorem}[section]
\newtheorem{lemma}[theorem]{Lemma}
\newtheorem{proposition}[theorem]{Proposition}
\newtheorem{corollary}[theorem]{Corollary}
\newtheorem*{proposition*}{Proposition}

\theoremstyle{definition}
\newtheorem{remark}[theorem]{Remark}
\newtheorem{definition}[theorem]{Definition}
\newtheorem{example}[theorem]{Example}
\newtheorem*{notation}{Notation}
\newtheorem*{acknowledgements}{Acknowledgements}

\numberwithin{equation}{section}

\newcommand{\A}{{\mathfrak A}}
\newcommand{\cA}{{\mathcal A}}
\newcommand{\fA}{{\mathfrak A}}

\newcommand{\fC}{{\mathfrak C}}

\newcommand{\cG}{{\mathcal G}}
\newcommand{\bG}{{\mathbf G}}
\newcommand{\fG}{{\mathfrak G}}

\newcommand{\bH}{{\mathbf{H}}}

\renewcommand{\L}{{\mathcal L}}
\newcommand{\cL}{{\mathcal L}}

\newcommand{\cM}{{\mathcal M}}
\newcommand{\fM}{{\mathfrak M}}

\newcommand{\cN}{{\mathcal N}}
\newcommand{\nr}{{\mathrm{nr}}}
\renewcommand{\O}{{\mathcal O}}

\newcommand{\brp}{{\{p\}}}

\newcommand{\Q}{{\mathbf Q}}
\newcommand{\bQ}{{\mathbf Q}}

\newcommand{\br}{{\mathbf r}}

\newcommand{\ti}[1]{{\tilde{#1}}}

\newcommand{\Up}{{\Upsilon}}

\newcommand{\Z}{{\mathbf Z}}
\newcommand{\bZ}{{\mathbf Z}}

\DeclareMathOperator{\Aut}{Aut}

\DeclareMathOperator{\Cl}{Cl}

\DeclareMathOperator{\Cores}{Cores}

\DeclareMathOperator{\ev}{ev}
\DeclareMathOperator{\Ext}{Ext}
\DeclareMathOperator{\Gal}{Gal}
\DeclareMathOperator{\GL}{GL}
\DeclareMathOperator{\Hom}{Hom}

\DeclareMathOperator{\Ker}{Ker}
\DeclareMathOperator{\Loc}{Loc}
\DeclareMathOperator{\Map}{Map}

\DeclareMathOperator{\ord}{ord}
\DeclareMathOperator{\Pic}{Pic}

\DeclareMathOperator{\red}{red}
\DeclareMathOperator{\Res}{Res}

\DeclareMathOperator{\Sp}{Sp}
\DeclareMathOperator{\Spec}{Spec}

\begin{document}
\title{Galois modules and $p$-adic representations}
\author{A. Agboola}
\date{Version of October 8, 2004}
\address{Department of Mathematics \\
         University of California \\
         Santa Barbara, CA 93106. }
\email{agboola@math.ucsb.edu}
\subjclass{11Gxx, 11Rxx}

\begin{abstract}
In this paper we develop a theory of class invariants associated to
$p$-adic representations of absolute Galois groups of number
fields. Our main tool for doing this involves a new way of describing
certain Selmer groups attached to $p$-adic representations in terms of
resolvends associated to torsors of finite group schemes. 
\end{abstract}

\maketitle


\section{Introduction} \label{S:intro}

In this paper we shall introduce and study invariants which measure
the Galois structure of certain torsors that are constructed via
$p$-adic Galois representations. We begin by describing the background
to the questions that we intend to discuss.

Let $Y$ be any scheme, and suppose that $G \to Y$ is a
finite, flat, commutative group scheme. Write $G^*$ for the Cartier
dual of $G$. Let $\tilde{G}^*$ denote the normalisation of $G^*$, and
let $i:\tilde{G}^* \to G^*$ be the natural map. Suppose that $\pi:X \to
Y$ is a $G$-torsor, and write $\pi_0:G \to Y$ for the trivial
$G$-torsor. Then $\O_X$ is an $\O_G$-comodule, and so it is also an
$\O_{G^*}$-module (see e.g. \cite{CEPT}). As an $\O_{G^*}$-module, the
structure sheaf $\O_X$ is locally free of rank one, and so it gives a
line bundle $\cM_{\pi}$ on $G^*$. Set
\begin{equation*}
\L_{\pi}:= \cM_{\pi} \otimes \cM_{\pi_0}^{-1}.
\end{equation*}
Then the maps
\begin{equation} \label{E:psi}
\psi: H^1(Y,G) \to \Pic(G^*),\qquad [\pi] \mapsto [\L_{\pi}];
\end{equation}
\begin{equation} \label{E:varphi}
\varphi: H^1(Y,G) \to \Pic(\widetilde{G}^*),\qquad [\pi] \mapsto
[i^*\L_{\pi}]
\end{equation}
are homomorphisms which are often referred to as `class invariant
homomorphisms'. 

The initial motivation for studying class invariant homomorphisms
arose from Galois module theory.  Let $F$ be a number field with ring
of integers $O_F$, and suppose that $Y = \Spec(O_F)$. Write $G^* =
\Spec(A)$, $G = \Spec(B)$, and $X= \Spec(C)$. Then the algebra $C$ is
a twisted form of $B$, and the homomorphisms $\psi$ and $\varphi$
measure the Galois module structure of this twisted form. The
homomorphism $\psi$ was first introduced by W. Waterhouse (see
\cite{W}), and was further developed in the context of Galois module
theory by M. Taylor (\cite{T1}). Taylor originally considered the case
in which $G$ is a torsion subgroup scheme of an abelian variety with
complex multiplication. The corresponding torsors are obtained by
dividing points in the Mordell-Weil groups of such abelian varieties,
and they are closely related to rings of integers of abelian
extensions of $F$. In \cite{ST}, it was shown that, for elliptic
curves with complex multiplication, the class invariant homomorphism
$\psi$ vanishes on the classes of torsors obtained by dividing torsion
points of order coprime to $6$. This implies the existence of Galois
generators for certain rings of integers of abelian extensions of
imaginary quadratic fields, and it may be viewed as an integral
version of the Kronecker Jugendtraum (see \cite{T2}, \cite{CNT}). This
vanishing result was extended to all elliptic curves in \cite{A3} and
\cite{P1}.

Since their introduction, class invariants of torsors obtained by
dividing points on abelian varieties have been studied in greater
generality by several authors. For example, suppose that $\mathcal{X}$
is a projective curve over $\Spec(\Z)$ which is equipped with a free
action of a finite group. In \cite{P2}, it is shown that the behaviour
of the equivariant projective Euler characteristic of
$\O_{\mathcal{X}}$ is partly governed by class invariants of torsors
arising from torsion points on the Jacobian of $\mathcal{X}$. In
\cite{AP}, an Arakelov (i.e. arithmetic) version of class invariants
of torsors coming from points on abelian varieties is
considered. There it is shown that in general such torsors are
completely determined by their arithmetic class invariants, and that
these invariants are related to Mazur-Tate heights on the abelian
variety (see \cite{MT}). Finally we mention that in \cite{A2},
\cite{A4}, and \cite{AT}, class invariants arising from points on
elliptic curves with complex multiplication are studied using Iwasawa
theory, and they are shown to be closely related to the $p$-adic
height pairing on the elliptic curve.

The main goal of this paper is to develop a theory of class invariants
for arbitrary $p$-adic representations, and to generalise a number of
results that up to now have only been known in certain cases involving
elliptic curves with complex multiplication.

We now describe the main results contained in this paper. Let $V$ be a
$d$-dimensional $\Q_p$-vector space. Let $F^c$ be an algebraic closure
of $F$, and write $\Omega_F:= \Gal(F^c/F)$. Suppose that
$\rho: \Omega_F \to \GL(V)$
is a continuous representation of $\Omega_F$ that is ramified at only
finitely many primes of $F$.
Set $V^*:= \Hom_{\Q_p}(V,\Q_p(1))$, and let
$\rho^*: \Omega_F \to \GL(V^*)$
be the corresponding representation of $\Omega_F$. Suppose that $T
\subseteq V$ is an $\Omega_F$-stable lattice, and write $T^*:=
\Hom_{\Z_p}(T,\Z_p(1))$. (Note that for each construction in this
paper that depends upon $T$, there is also a corresponding
construction that depends upon $T^*$; this will not always be
explicitly stated.)

For each positive integer $n$, we may define finite group schemes
$G_n$ and $G_n^*$ over
$\Spec(F)$ by
\begin{equation*}
G_n(F^c) =\Gamma_n:= p^{-n}T/T; \qquad G_n^*(F^c) = \Gamma_n^*:=
p^{-n}T^*/T^*.
\end{equation*}
Then $G_n^*$ is the Cartier dual of $G_n$, and we may write
$G_n^*=\Spec(A_n)$ for some Hopf algebra $A_n$ over $F$. Let $\fA_n$
be any $O_F$-algebra such that $\fA_n \otimes_{O_F}F = A_n$, and write
$\Sp(\fA_n):= \Spec(\fA_n)$.

By using a description of $H^1(F,G_n)$ which arises via studying the
Galois structure of $G_n$-torsors in terms of $A_n$, we give a new way
of imposing local conditions on cohomology classes in terms of the
algebra $\fA_n$. (Roughly speaking, if $\pi \in H^1(F,G_n)$, then we
use $\fA_n$ to impose local conditions on line bundle $\L_{\pi}$
associated to $\pi$.) This yields a certain Selmer group in
$H^1(F,G_n)$ which we denote by $H^{1}_{\fA_n}(F,G_n)$. Suppose that
$\pi:X \to \Spec(F)$ is any $G_n$-torsor whose isomorphism class lies
in $H^{1}_{\fA_n}(F,G_n)$. We shall explain how to use the methods of
\cite{By}, \cite{Mc} and \cite{W} to construct a natural homomorphism
\begin{equation} \label{E:geninv}
\phi_{\fA_n}: H^{1}_{\fA_n}(F,G_n) \to \Pic(\Sp(\fA_n)).
\end{equation}
This generalises the class invariant homomorphisms \eqref{E:psi} and
\eqref{E:varphi} above.  For suppose that $G_n$ is the generic fibre of a
finite, flat group scheme $\cG_n$ over $\Spec(O_F)$. If we choose
$\fA_n$ to be the $O_F$-Hopf algebra representing the Cartier dual
$\cG_n^*$ of $\cG_n$, then $H^{1}_{\fA_n}(F,G_n) = H^1(\Spec(O_F),
\cG_n)$, and $\phi_{\fA_n}$ is the same as the homomorphism
\eqref{E:psi} in this case. If on the other hand we take $\fA_n$ to be
the maximal $O_F$-order $\fM_n$ in $A_n$, then $\Sp(\fA_n)$ is equal
to the normalisation $\tilde{\cG}_n^*$ of $\cG^*_n$. In this case,
$H^1(\Spec(O_F), \cG_n)$ is contained in $H^{1}_{\fA_n}(F,G_n)$, and
the restriction of $\phi_{\fA_n}$ to $H^1(\Spec(O_F), \cG_n)$ is the
homomorphism \eqref{E:varphi}. (See Example \ref{e:flat} below.)

In this paper we shall mainly be concerned with the cases
$$
\fA_n = \fM_n,\qquad \fA_n = \fM_n \otimes_{O_F} O_F[1/p]:=
\fM_{n}^{\brp}.
$$
For each finite place $v$ of $F$, let $F^\nr_v$ denote the maximal
unramified extension of $F_v$ in a fixed algebraic closure of
$F_v$. If $v \nmid p$, then define
\begin{equation} \label{E:f unram}
H^1_f(F_v,T):= \Ker \left[ H^1(G_v,T) \to H^1(F^\nr_v,T)
\right]. 
\end{equation}
Following \cite[\S3.1.4]{PR1}, we set
\begin{equation*}
H^1_{f,\{p\}}(F,T) = \Ker \left[ H^1(F,T) \to \oplus_{v \nmid p}
\frac{H^1(F_v,T)}{H^1_f(F_v,T)} \right].
\end{equation*}
It may be shown that (see Remark \ref{R:unram} below)
$$
H^1_{f,\{p\}}(F,T) \subseteq \varprojlim H^{1}_{\fM_{n}^{\brp}}(F,G_n)
$$
Here the inverse limit is taken with respect to the maps induced by
the `multiplication by $p$' maps $G_{n+1} \to G_{n}$, and we view
$\varprojlim H^{1}_{\fM_{n}^{\brp}}(F,G_n)$ as being a subgroup of
$H^1(F,T)$ via the canonical isomorphism $\varprojlim H^1(F,G_n)
\simeq H^1(F,T)$. Set $H^1_u(F,T):= \varprojlim H^{1}_{\fM_n}(F,G_n)$.

The natural inclusions $G_n^* \to G_{n+1}^{*}$ induce pullback
homomorphisms
$$
\Pic(\Sp(\fM_{n+1})) \to \Pic(\Sp(\fM_n)),\qquad
\Pic(\Sp(\fM_{n+1}^{\brp})) \to \Pic(\Sp(\fM_{n}^{\brp})).
$$
We shall show that we may take inverse limits in \eqref{E:geninv} to
obtain homomorphisms
$$
\Phi_F: H^1_u(F,T) \to \varprojlim \Pic(\Sp(\fM_n)),\qquad
\Phi_{F}^{\brp}: H^{1}_{f,\brp}(F,T) \to \varprojlim
\Pic(\Sp(\fM_{n}^{\brp})).
$$

Our first result shows that the homomorphism $\Phi_{F}^{\brp}$ is
closely related to a $p$-adic height pairing associated to $T$. In
order to describe why this is so, we have to introduce some further
notation.

Let $C_n/F$ denote the $n$-th layer of the cyclotomic
$\bZ_p$-extension of $F$, and set
$$
\fG_F(T):= \{ x \in H^{1}_{f,\brp}(F,T) \mid Mx \in \cap_n
\Cores_{C_n/F}(H^{1}_{f,\brp}(C_n,T)) \, \text{for some integer
$M>0$} \}.
$$
When $T$ is the $p$-adic Tate module of an elliptic curve, $\fG_F(T)$
is the same as the canonical subgroup that was defined by R. Greenberg
in \cite[p.131--132]{G}, and further studied by A. Plater in \cite{Pl}
and \cite{Pl2} (see also \cite{JJ}).

Let
\begin{equation*}
\Loc_{F,T^*}: H^1_{f,\{p\}}(F,T^*) \to \prod_{v \mid p} H^1(F_v,T^*)
\end{equation*}
denote the natural localisation map. In \cite[Section 3.1.4]{PR1}, Perrin-Riou
constructs a $p$-adic height pairing
\begin{equation*}
B_F: H^1_{f,\{p\}}(F,T) \times \Ker(\Loc_{F,T^*}) \to \Q_p,
\end{equation*}
and she shows that the group $\fG_F(T)$ lies in the left-hand kernel of
$B_F$. Write
\begin{equation} \label{E:height}
\langle \langle\,,\,\rangle \rangle:
\frac{H^1_{f,\{p\}}(F,T)}{\fG_F(T)} \times \Ker(\Loc_{F,T^*}) \to \Q_p
\end{equation}
for the pairing induced by $B_F$. It is conjectured that $\langle
\langle\,,\,\rangle \rangle$ is always non-degenerate modulo
torsion. If this conjecture is true, then it implies that $\fG_F(T)$
has a natural characterisation in terms of $p$-adic height pairings
attached to $T$.  The following result shows that this conjecture
implies that $\fG_F(T)$ also has a natural characterisation in terms of
Galois module structure.

\begin{theorem} \label{T:intro 1}
Suppose that the pairing $\langle \langle\,,\,\rangle \rangle$ is
non-degenerate modulo torsion. Then $x \in \fG_F(T)$ if and only if
$\Phi_{F}^{\brp}(x)$ is of finite order.
\end{theorem}

Let $S$ be a finite set of places of $F$ containing the places lying
over $p$, the places at which $\rho$ is ramified, and the set of
infinite places. Let $F^S/F$ denote the maximal unramified outside $S$
extension of $F$. A conjecture of Greenberg asserts that the group 
$H^2(F^S/C_\infty, V^*/T^*)$ always vanshes. This may be viewed as an
analogue of a weak form of Leopoldt's conjecture for the Galois
representation $V^*$.

\begin{corollary} \label{C:intro}
Suppose that the pairing $\langle \langle\,,\,\rangle \rangle$ is
non-degenerate modulo torsion, and that 
\begin{equation} \label{E:ralph}
H^2(F^S/C_\infty, V^*/T^*) = 0.
\end{equation}
Then the restriction of $\Phi_{F}^{\brp}$ to $\Ker(\Loc_{F,T})$ has
finite kernel.
\end{corollary}

\begin{proof} It follows from \cite[Remark at the end of
\S3.4.3]{PR1} that if \eqref{E:ralph} holds, then $\fG_F(T) \cap
\Ker(\Loc_{F,T})$ is finite. The result now follows from Theorem
\ref{T:intro 1}.
\end{proof}

We now turn to the homomorphism $\Phi_F$.

For each integer $n$, the action of $\Omega_F$ on $\Gamma_n^*$ yields
a representation $\rho_n^*: \Omega_F \to \Aut(\Gamma_n^*)$. Write
$F_n^*$ for the fixed field of $\rho_n^*$; then $F^*_\infty:= \cup_n
F_n^*$ is the extension of $F$ cut out by $\rho^*$. Set
\begin{equation*}
\fC_F(T):= \{ x \in H^1_u(F,T) \mid x \in \cap_n
\Cores_{F_{n}^{*}/F}(H^1_u(F_n^*,T)) \, \text{for some integer $M>0$}
\}.
\end{equation*}

\begin{theorem} \label{T:intro 2}
Suppose that $x \in \fC_F(T)$. Then $\Phi_F(x)$ is of finite order.
\end{theorem}

\begin{remark} \label{R:intro}
Whether or not the the converse of Theorem \ref{T:intro 2} holds in
general is an open question, and it appears to be a very delicate
problem. If $T$ has $\bZ_p$-rank one, then it may be shown that
$\fG_F(T) = \fC_F(T)$, and so the converse to Theorem \ref{T:intro 2}
holds in this case. \qed
\end{remark}


\begin{acknowledgements} This paper owes a great deal to ideas first
introduced in \cite{AP}, and I am extremely grateful to G. Pappas for
many interesting and helpful conversations. I would like to thank
B. Conrad, R. Greenberg and S. Howson for useful discussions, and an
anonymous referee for pointing some errors in an earlier version of
this paper. I would also like to thank the Mathematics Department of
Harvard University for their hospitality while a part of this work was
carried out. This research was partially supported by NSF grants.
\end{acknowledgements}

\begin{notation}
In this paper, all modules are left modules unless explicitly stated
otherwise.

For any field $L$, we write $L^c$ for an algebraic closure of $L$, and
we set $\Omega_L:= \Gal(L^c/L)$. If $L$ is either a number field or a
local field, then we write $O_L$ for its ring of integers.

If $L$ is a number field and $v$ is a finite place of $L$, then we
write $L_v$ for the local completion of $L$ at $v$. We fix
an algebraic closure $L_v^c$ of $L_v$ and we identify $\Omega_{L_v}$
with a subgroup of $\Omega_L$. If $P$ is any $O_L$-module, then we
shall usually write write $P_v:= P \otimes_{O_L} O_{L_v}$. 

For any $\Z$-module $Q$, we set $\Check{Q}:= \varprojlim_n Q/p^n Q$.

If $R$ and $S$ are rings with $R \subseteq S$, and if $A$ is any
$R$-algebra, then we often write $A_S$ for $A \otimes_R S$. We shall
also often write $\Sp(A)$ for $\Spec(A)$.
\end{notation}


\section{Resolvends and cohomology groups} \label{S:cohomology}


Let $R$ be a Dedekind domain with field of fractions $K$, and write
$R^c$ for the integral closure of $R$ in $K^c$. (In what follows we
always allow the possibility that $R=K$ and $R^c=K^c$.) Assume that
$K$ is of characteristic zero. In this section, we shall take $Y =
\Spec(R)$, and we shall explain how the cohomology group $H^1(Y,G)$
may be described in terms of the Hopf algebra $A$ representing
$G^*$. Set $\Gamma:= G(R^c)$ and $\Gamma^*:= G^*(R^c)$.

Recall that there is a canonical isomorphism
\begin{equation*}
H^1(Y,G) \simeq \Ext^1(G^*,\bG_m)
\end{equation*}
(see \cite{W}, \cite[\'expos\'e VII]{SGA7}, or \cite{P2}). This
implies that given any $G$-torsor $\pi:X \to Y$, we can associate to it
a canonical commutative extension
\begin{equation*}
1 \to \bG_m \to G(\pi) \to G^* \to 1.
\end{equation*}
The scheme $G(\pi)$ is a $\bG_m$-torsor over $G^*$, and its associated
$G^*$-line bundle is equal to $\L_{\pi}$. (This construction is
explained in detail by Waterhouse in \cite{W}.)

Over $\Spec(R^c)$, the $G$-torsors $\pi_0$ and $\pi$ become
isomorphic, i.e. there is an isomorphism $X \otimes_R R^c \simeq G
\otimes_R R^c$ of schemes with $G$-action. (This isomorphism is not
unique: it is only well-defined up to the action of an element of
$G(R^c)$.) Hence, via the functoriality of Waterhouse's construction in
\cite{W}, we obtain an isomorphism
\begin{equation*}
\xi_{\pi}: \L_{\pi} \otimes_R R^c \xrightarrow{\sim} A_{R^c}.
\end{equation*}
We shall refer to $\xi_{\pi}$ as a \textit{splitting isomorphism} for
$\pi$.

Now suppose that $\psi(\pi)=0$. Then $\L_{\pi}$ is a free $A$-module,
and so we may choose a trivialisation $s_{\pi}:A \xrightarrow{\sim}
\L_{\pi}$. Consider the composition
\begin{equation*}
A_{R^c} \xrightarrow{s_{\pi} \otimes_R R^c} \L_{\pi} \otimes_R R^c
\xrightarrow{\xi_{\pi}} A_{R^c}.
\end{equation*}
This is an isomorphism of $A_{R^c}$-modules, and so it is just
multiplication by an element $\br(s_{\pi})$ of $A_{R^c}^{\times}$. We
refer to $\br(s_{\pi})$ as a {\it resolvend} of $s_{\pi}$ or as a {\it
resolvend associated to} $\pi$. (This terminology is due to
L. McCulloh, \cite{Mc}.) Note that $\br(s_{\pi})$ depends upon the
choice of $\xi_{\pi}$ as well as upon $s_{\pi}$.

\begin{definition} Let $A$ and $R$ be as above. Define
\begin{align*}
\bH(A) &:= \{ \alpha \in A_{R^c}^{\times}\,|\, \alpha^{\omega} = g_{\omega}
\alpha \,\, \text{for all $\omega \in \Omega_K$, where $g_{\omega} \in
\Gamma$}\}, \\
H(A)&:= \frac{\bH(A)}{\Gamma \cdot A^{\times}}. 
\end{align*}
\qed
\end{definition}

If $\omega \in \Omega_K$, then $\xi_{\pi}^{\omega} = g_{\omega}
\xi_{\pi}$, where $g_{\omega} \in \Gamma$. Since $s_{\pi}^{\omega} =
s_{\pi}$, we deduce that $\br(s_{\pi})^{\omega} = g_{\omega}\br(s_{\pi})$,
that is, $\br(s_{\pi})^{\omega} \in \bH(A)$. It is easy to see that
changing $s_{\pi}$ alters $\br(s_{\pi})$ via multiplication by an
element of $A^{\times}$, while changing $\xi_{\pi}$ alters
$\br(s_{\pi})$ via multiplication by an element of $\Gamma$. Hence the
image $r(\pi)$ of $\br(s_{\pi})$ in $H(A)$ depends only upon
the isomorphism class of the torsor $\pi$.

The following result, in the case in which $G$ is a constant group
scheme, is equivalent to certain results of L. McCulloh (see
\cite[Sections 1 and 2]{Mc}; note, however that McCulloh formulates
his results in a rather different way from that described
here). McCulloh's methods were generalised by N. Byott to the case of
arbitrary $G$ (see \cite[Lemma 1.11 and Sections 2 and 3]{By}) using
techniques from the theory of Hopf algebras. We give a different
approach to these ideas.

\begin{theorem} \label{T:Leon}
Let $G$ be a finite, flat commutative group scheme over $\Spec(R)$,
and let $G^* = \Spec(A)$ be the Cartier dual of $G$. Then the map
\begin{equation*}
\Up_R: \Ker(\psi) \to H(A);\qquad
[\pi] \mapsto r(\pi)
\end{equation*}
is an isomorphism.
\end{theorem}

\begin{proof} 
We first show that $\Up_R$ is a homomorphism. Suppose that $\pi:X \to
\Spec(R)$ and $\pi':X' \to \Spec(R)$ are $G$-torsors satisfying
$\psi(\pi) = \psi(\pi') = 0$. Let $s_{\pi}: A \xrightarrow{\sim}
\L_{\pi}$ and $s_{\pi'}: A \xrightarrow{\sim} \L_{\pi'}$ be
trivialisations of $\L_{\pi}$ and $\L_{\pi'}$ respectively. Write
$\pi'':= \pi \cdot \pi'$. Then it follows via the functoriality of
Waterhouse's construction in \cite{W} that there is a natural
isomorphism $\L_{\pi''} \simeq \L_{\pi} \otimes \L_{\pi'}$. Thus, if
we set
\begin{equation*}
s_{\pi''} = s_{\pi} \otimes s_{\pi'}: A \simeq A \otimes_A A
\xrightarrow{\sim} \L_{\pi} \otimes_A \L_{\pi'},
\end{equation*}
then we see that $\br(s_{\pi''}) = \br(s_{\pi}) \br(s_{\pi'})$. This
implies that $r(\pi'') = r(\pi) r(\pi')$, as required.

We now show that $\Up_R$ is surjective. For any scheme $S \to
\Spec(R)$, write $\Map_S(G^*, \bG_m)$ for the set of scheme-theoretic
morphisms $G^* \otimes_R S \to \bG_m \otimes_R S$. Since $G^*$ is
affine, the functor $S \mapsto \Map_S(G^*, \bG_m)$ is representable by
the affine group scheme $\bG_{m,G^*}:= \bG_m \otimes_R G^* \to
\Spec(R)$. The group scheme $G= \Hom(G^*,\bG_m)$ is a closed subgroup
scheme of $\bG_{m,G^*}$.

Suppose that $\alpha \in \bH(A)$. Then we may view $\alpha$ as being a
$\Spec(R^c)$-valued point of $\bG_{m,G^*}$. Let
\begin{equation}
G_{\alpha}:= \alpha \cdot [G \otimes_R R^c]
\end{equation}
denote the `translation-by-$\alpha$' in $\bG_{m,G^*} \otimes_R R^c$ of $G
\otimes_R R^c$. Then $G \otimes_R R^c$ acts on $G_{\alpha}$ via
translation. Furthermore, translation by $\alpha$ induces an
isomorphism
\begin{equation}
\Xi_{\alpha}: G \otimes_R R^c \xrightarrow{\sim} G_{\alpha}
\end{equation}
of schemes with $G \otimes_R R^c$ action.

We now claim that $G_{\alpha}$ descends to $\Spec(R)$, i.e. that there
is a scheme $\pi_{\alpha}:Z_{\alpha} \to \Spec(R)$ defined over $R$
which is such that $G_{\alpha} = Z_{\alpha} \otimes_R R^c$. Since $R$
is a Dedekind domain, and $G_{\alpha}$ is flat over $\Spec(R^c)$, it
suffices to check that the generic fibre $G_{\alpha/{K^c}}$ of
$G_{\alpha}$ descends to a scheme over $\Spec(K)$. This in turn
follows via Galois descent, and may be seen as follows. We first note
that the isomorphism $\Xi_{\alpha}$ induces a bijection $\Gamma \to
G_{\alpha}(K^c)$ of sets. Define $z_{\alpha}: \Omega_K \to \Gamma$ by
$z(\omega) = \alpha^{\omega} \alpha^{-1}$; thus $z_{\alpha}$ is the
$\Gamma$-valued cocycle of $\Omega_K$ associated to $\alpha$. Then it
is easy to check that the action of $\Omega_K$ on $G_{\alpha}(K^c)$ is
given by
\begin{equation*}
\Xi(g)^{\omega} = z_{\alpha}(\omega) \Xi(g^{\omega})
\end{equation*}
for all $g \in \Gamma$ and $\omega \in \Omega_K$. This implies that
$G_{\alpha/{K^c}}$ descends to $Z_{\alpha/K}$ over $\Spec(K)$. A
similar argument also shows that $\pi_{\alpha}:Z_{\alpha} \to
\Spec(R)$ is a $G$-torsor over $\Spec(R)$. 

We shall now show that $\psi(\pi_{\alpha}) = 0$. Let
\begin{equation*}
\xi_{\pi_{\alpha}}: \L_{\pi_{\alpha}} \otimes_R R^c \xrightarrow{\sim}
A_{R^c}
\end{equation*}
denote the splitting isomorphism of $\pi_{\alpha}$ induced by
$\Xi_{\alpha}$. Define an isomorphism
\begin{equation*}
\sigma_{\alpha}: A_{R^c} \xrightarrow{\sim} \L_{\pi_\alpha} \otimes_R R^c
\end{equation*}
by $\sigma_{\alpha}(a) = \xi_{\pi_{\alpha}}^{-1}(\alpha a) = \alpha
\xi_{\pi_{\alpha}}^{-1}(a)$ for all $a \in A_{R^c}$. In order to show that
$\psi(\pi_{\alpha}) = 0$, it suffices to show that $\sigma_{\alpha}$
descends to an isomorphism $\sigma_{\alpha}':A \xrightarrow{\sim}
\L_{\pi_\alpha}$ over $R$. This will in turn follow if we show that
\begin{equation*}
\sigma_{\alpha}^{\omega}(a) = \sigma_{\alpha}(a)
\end{equation*}
for all $\omega \in \Omega_K$ and all $\alpha \in A_{R^c}$. To check
this last equality, we simply observe that
\begin{align*}
\sigma_{\alpha}^{\omega}(a) = \omega
[\sigma_{\alpha}(a^{\omega^{-1}})] 
&= \omega [ \alpha \xi_{\pi_{\alpha}}^{-1}( a^{\omega^{-1}})] \\
&= \omega [ \alpha z_{\alpha}(\omega^{-1})
\xi_{\pi_{\alpha}}^{-1}(a)^{\omega^{-1}}]  \\
&= \omega [ \alpha^{\omega^{-1}}
\xi_{\pi_{\alpha}}^{-1}(a)^{\omega^{-1}}] \\
&= \alpha \xi_{\pi_{\alpha}}^{-1}(a) \\
&= \sigma_{\alpha}(a).
\end{align*} 
Hence $\psi(\pi_{\alpha}) = 0$ as asserted.

To complete the proof of the surjectivity of $\Up_R$, we note that it
follows from the definition of $\sigma_{\alpha}$ that we have
$\br(\sigma_{\alpha}') = \alpha$. Hence $r(\pi_{\alpha}) = [\alpha]
\in H(A)$, and so $\Up_R$ is surjective as claimed.

We now show that $\Up_R$ is injective. Suppose that $\alpha, \beta \in
\bH(A)$ with $[\alpha] = [\beta] \in H(A)$. Then it is easy
to check that the isomorphism
\begin{equation*}
\Xi_{\beta} \circ \Xi_{\alpha}^{-1}: G_{\alpha} \xrightarrow{\sim}
G_{\beta}
\end{equation*}
induces an isomorphism $G_{\alpha}(K^c) \xrightarrow{\sim} G_{\beta}(K^c)$
of $\Omega_K$-modules. This implies that the $G$-torsors
$\pi_{\alpha}:Z_{\alpha} \to \Spec(R)$ and $\pi_{\beta}: Z_{\beta} \to
\Spec(R)$ are isomorphic. This completes the proof of the theorem.
\end{proof}

\begin{remark} \label{R:cocycle}
Suppose that $R=K$. Then it is not hard to check that
(using the notation established in the proof of Theorem \ref{T:Leon})
the map $\Omega_K \to \Gamma$ defined by $\omega \mapsto
\br(s_{\pi})^{\omega} \br(s_{\pi})^{-1}$ is an $\Omega_K$-cocycle
representing $[\pi] \in H^1(K,G)$. \qed
\end{remark}

If $R$ is a local ring, then $\Pic(G^*)=0$, and so $\Ker(\psi) =
H^1(R, G)$. The following result is a direct corollary of
Theorem \ref{T:Leon}. It gives a description of the flat cohomology of
$G$ over $\Spec(R)$ in terms of resolvends. (Recall that $N$ denotes
the exponent of $G$.)

\begin{corollary} \label{C:Leon} Suppose that $R$ is a local ring.

(a) There is an isomorphism
\begin{equation*}
\Up_R:H^1(R,G) \xrightarrow{\sim} H(A).
\end{equation*}

(b) The map $[\pi] \mapsto \br(s_{\pi})^N$ induces a homomorphism
\begin{equation*}
\eta_R: H^1(R,G) \to \frac{A^{\times}}{(A^{\times})^N}. 
\end{equation*}
\qed
\end{corollary}

\begin{remark} \label{R:wedd}
Suppose that $R=K$, and for each $\gamma^* \in \Gamma^*$, write
$K[\gamma^*]$ for the smallest extension of $K$ whose absolute Galois
group fixes $\gamma^*$. Let $\Gamma^* \backslash \Omega_K$ denote a
set of representatives of $\Omega_K$-orbits of $\Gamma^*$. Then, via
an argument virtually identical to that given in \cite[Lemma 3.3]{A2},
it may be shown that the Wedderburn decomposition of the $K$-algebra
$A$ is given by
\begin{equation} \label{E:wedd}
A \simeq (K^c \Gamma)^{\Omega_K}
\simeq \prod_{\gamma^* \in \Gamma^* \backslash \Omega_K}
K[\gamma^*].
\end{equation}
\qed
\end{remark}

\begin{proposition} \label{P:etaprim}
Suppose that $R=K$, and that $G^*$ is a constant group scheme over
$\Spec(K)$. Then $A \simeq \Map(\Gamma^*,K)$, and the map $\eta_K$ of
Corollary \ref{C:Leon}(b) induces an isomorphism
\begin{equation*}
\eta_K:H^1(K,G) \xrightarrow{\sim}
\Hom(\Gamma^*,K^{\times}/(K^{\times})^N) \subset
\Map(\Gamma^*,K^{\times}/(K^{\times})^N) \simeq
A^{\times}/(A^{\times})^N.
\end{equation*}
\end{proposition}

\begin{proof} See \cite[Corollary 3.4]{AP}.
\end{proof}
 
\begin{proposition} \label{P:incres}
Suppose that $R=K$, and let $L$ be a (possibly infinite) extension of
$K$. Then the following diagram is commutative:
\begin{equation} \label{E:incres}
\begin{CD}
H^1(K,G) @>{\Up_K}>> H(A) \\
@V{\Res}VV           @VVV     \\
H^1(L,G)  @>{\Up_L}>> H(A_L).
\end{CD}
\end{equation}
Here the left-hand vertical arrow is the restriction map on
cohomology, and the right-hand vertical arrow is the homomorphism
induced by the inclusion map $i: \bH(A) \to \bH(A_L)$.
\end{proposition}

\begin{proof}
Let $\pi:X \to \Spec(K)$ be any $G$-torsor, and let $s:A
\xrightarrow{\sim} \cL_\pi$ be any trivialisation of $\cL_\pi$. Then
it follows via a straightforward computation that the
$\Omega_L$-cocycle associated to $i(\br(s))$ is equal to the
restriction of the $\Omega_K$-cocycle associated to $\br(s)$
(cf. Remark \ref{R:cocycle}).
\end{proof}

\begin{remark} \label{R:kereta}
Suppose that $R=K$, and that $\pi \in \Ker(\eta_F)$. Let $\br(s_\pi)
\in \bH(A)$ be any resolvend associated to $\pi$. Then $\br(s_\pi)^N =
\alpha^N \in A^{\times N}$, and so $\br(\alpha^{-1}s_\pi)^N =
1$. Hence $\br(\alpha^{-1} s_\pi) \in A_{F(\mu_N)}^{\times}$, and so
Proposition \ref{P:incres} implies that $\pi$ lies in the kernel of
the restriction map 
$$
\Res_{K/K(\mu_N)}: H^1(K,G) \to H^1(K(\mu_N),G).
$$

Conversely, if $\pi \in \Res_{K/K(\mu_N)}$, then for a suitable choice
of $s_\pi$, we have $\br(s_\pi) \in A_{F(\mu_N)}^{\times}$, and
$\br(s_\pi)^N \in A^{\times}$. Since $A$ is a semisimple $F$-algebra,
this implies (via considering the Wedderburn decomposition
\eqref{E:wedd} of $A$) that $\br(s_\pi)^N \in A^{\times N}$, and so
$\pi \in \Ker(\eta_F)$.  \qed
\end{remark}

Suppose now that $L$ is a finite Galois extension of $K$ with
$[L:K]=n$, say. Let $\omega_1, \ldots, \omega_n$ be a transversal of
$\Omega_L$ in $\Omega_K$. Then we have a norm homomorphism
\begin{equation} \label{E:norm}
\cN_{L/K}: A_{L^c} \to A_{K^c};\qquad a \mapsto \prod_{i=1}^{n}
a^{\omega_i}.
\end{equation}
This induces homomorphisms (which we denote by the same symbol)
\begin{equation*}
\cN_{L/K}: \bH(A_L) \to \bH(A_K),\qquad \cN_{L/K}: H(A_L) \to H(A_K).
\end{equation*}

\begin{theorem} \label{T:traceres}
The following diagram is commutative:
\begin{equation} \label{E:traceres}
\begin{CD}
H^1(L,G) @>{\Up_L}>> H(A_L) \\
@V{\Cores_{L/K}}VV           @VV{\cN_{L/K}}V     \\
H^1(K,G)  @>{\Up_K}>> H(A_K),
\end{CD}
\end{equation}
where the left-hand vertical arrow is the corestriction map on
cohomology.
\end{theorem}

\begin{proof} Let $\pi_L: X_L \to \Spec(L)$ be any $G$-torsor and let
$s_{\pi_L}: A_L \xrightarrow{\sim} \L_{\pi_L}$ be 
any trivialisation of $\L_{\pi_L}$. Then it follows via a
straightforward computation that the $\Omega_K$-cocycle associated to
$\cN_{L/K}(\br(s_{\pi_L}))$ is equal to the corestriction of the
$\Omega_L$-cocycle associated to $\br(s_{\pi_L})$ (cf. Remark
\ref{R:cocycle}).
\end{proof}
\smallskip

Let $(G_n)_{n \geq 1}$ be a $p$-divisible group over $\Spec(R)$. For
each $n$, set $G_n^* = \Spec(A_n)$, and set $\Gamma_n:=G_n(R^c)$,
$\Gamma_n^*:=G_n^*(R^c)$. We write $p_n:=[p]: G_n \to G_{n-1}$ for the
multiplication-by-$p$ map, and we use the same symbol to denote the
induced map $H^1(\Spec(R),G_n) \to H^1(\Spec(R),G_{n-1})$.

The map $p_n$ induces a dual inclusion map $p_n^D:G_{n-1}^* \to
G_n^*$, and we may identify $A_{n-1}$ with the pullback $(p_n^D)^*A_n$
of $A_n$ via $p_n^D$. Thus (via pullback) $p_n^D$ induces a
homomorphism $q_n:A_n \to A_{n-1}$ which extends to a homomorphism
(which we denote by the same symbol) $A_{n,K^c} \to A_{n-1,K^c}$. It
is easy to check that $q_n(\bH(A_n)) \subseteq \bH(A_{n-1})$, and that
$q_n(\Gamma_n \cdot A_n^{\times}) \subseteq \Gamma_{n-1} \cdot
A_{n-1}^{\times}$.

\begin{theorem} \label{T:pdiv}
Suppose that $R$ is a local ring. Then the following
diagram is commutative:
\begin{equation} \label{E:pdiv}
\begin{CD}
H^1(\Spec(R),G_n) @>{\Up_{R,n}}>> H(A_n) \\
@V{p_n}VV           @VV{q_n}V     \\
H^1(\Spec(R),G_{n-1})  @>{\Up_{R,n-1}}>> H(A_{n-1}). 
\end{CD}
\end{equation}
\end{theorem}

\begin{proof} Suppose that $\pi_n:X_n \to \Spec(R)$ is any
$G_n$-torsor, and let $s_{\pi_n}: A_n \xrightarrow{\sim} \L_{\pi_n}$
be any trivialisation of $\L_{\pi_n}$. Set $\pi_{n-1}:=
p_n(\pi_n)$. Then it follows via the functoriality of Waterhouse's
construction in \cite{W} that there is a natural identification
$\L_{\pi_{n-1}} \simeq (p_n^D)^* \L_{\pi_n}$. Consider the
trivialisation $s_{\pi_{n-1}}:= (p_n^D)^*s_{\pi_n}: A_{n-1}
\xrightarrow{\sim} \L_{\pi_{n-1}}$ of $\L_{\pi_{n-1}}$ obtained by
pulling back $s_{\pi_n}$ along $p_n^D$.  We have
\begin{equation*}
\Up_{R,n-1}(\pi_{n-1}) = [\br(s_{\pi_{n-1}})] = [\br((p_n^D)^* s_{\pi_n})] =
[q_n(\br(s_{\pi_n})].
\end{equation*}
This establishes the result.
\end{proof}

Suppose now that $R=K$. Fix a positive integer $n$, and assume that
$G_n^*$ is a constant group scheme. Then we have
\begin{equation*}
A_{n}^{\times}/(A_{n}^{\times})^{p^n} \simeq \Map(\Gamma_n^*,
K^{\times}/(K^{\times})^{p^n}).
\end{equation*}
For each element $P: \Spec(K) \to G_n^*$ in $\Gamma_n^*$, write
$\chi_P:G_n \to \mu_{p^n}$ for the corresponding character of
$G_n$. Then $\chi_P$ induces a homomorphism (which we denote by the
same symbol):
\begin{equation*}
\chi_P: H^1(K,G_n) \to H^1(K,\mu_{p^n});\qquad [\pi] \mapsto
[\pi(\chi_P)].
\end{equation*}
Write
\begin{equation*}
\ev_P: A_{n}^{\times}/(A_{n}^{\times})^{p^n} \simeq \Map(\Gamma_n^*,
K^{\times}/(K^{\times})^{p^n}) \to K^{\times}/(K^{\times})^{p^n}
\end{equation*}
for the map $a \mapsto a(P)$ given by `evaluation at $P$'. The
following result shows how to describe the map $\eta_K$ of Corollary
\ref{C:Leon} in terms of Kummer theory.

\begin{proposition} \label{P:kummer key}
Let the hypotheses and notation be as above. Then the following
diagram is commutative:
\begin{equation}
\begin{CD}
H^1(K,G_n) @>{\chi_P}>> H^1(K,\mu_{p^n}) \\
@V{\eta_K}VV                    @AA{\mathrm{Kummer}}A     \\
A_n^{\times}/(A_n^{\times})^{p^n}     @>{\ev_P}>>
K^{\times}/(K^{\times})^{p^n}. 
\end{CD}
\end{equation}
(Here the right-hand vertical arrow is the natural isomorphism
afforded by Kummer theory.)
\end{proposition}

\begin{proof} See \cite[Proposition 3.2]{AP}.
\end{proof}

\begin{corollary} \label{C:pdiv compatible}

Let the hypotheses and notation be as above. For each
integer $n$, let 
\begin{equation*}
r_n: \Hom(\Gamma_n^*, K^{\times}/(K^{\times})^{p^n}) \to
\Hom(\Gamma_{n-1}^{*}, K^{\times}/(K^{\times})^{p^{n-1}})
\end{equation*}
be the homomorphism given by $f \mapsto f|_{\Gamma_{n-1}^{*}}$. Then the
following diagram commutes:
\begin{equation*}
\begin{CD}
H^1(K,G_n) @>{\eta_K}>> \Hom(\Gamma_n^*,K^{\times}/(K^{\times})^{p^n}) \\
@V{p_n}VV                     @VV{r_n}V              \\
H^1(K,G_{n-1}) @>{\eta_K}>>
\Hom(\Gamma_{n-1}^*,K^{\times}/(K^{\times})^{p^{n-1}}). 
\end{CD}
\end{equation*}
\end{corollary}

\begin{proof}
Suppose that $P \in \Gamma^{*}_{n-1}$. Then, almost by definition, the
following diagram commutes:
\begin{equation} \label{E:pdiv compatible}
\begin{CD}
H^1(K,G_n) @>{\chi_P}>> H^1(K,\mu_{p^n}) \simeq
K^{\times}/(K^{\times})^{p^n} \\ 
@V{p_n}VV                        @VV{\red}V \\
H^1(K,G_{n-1}) @>{\chi_P}>> H^1(K,\mu_{p^{n-1}}) \simeq
K^{\times}/(K^{\times})^{p^{n-1}}.
\end{CD}
\end{equation}

(Here the right-hand vertical arrow denotes the natural reduction
map.) The result now follows from Propositions \ref{P:kummer key} and
\ref{P:etaprim}.
\end{proof}


\section{Selmer conditions and Galois structure} \label{S:selmer}


In this section we shall apply the results of \S\ref{S:cohomology} to
explain how resolvends may be used to impose local conditions on
$G$-torsors. This enables us to define certain Selmer groups. We then
show that there are natural homomorphisms from these Selmer groups
into suitable locally free classgroups. These generalise the class
invariant homomorphisms described at the begining of the introduction
to this paper.
\smallskip

In what follows, $F$ will denote either a number field or a local
field (depending upon the context), with ring of integers $O_F$. We
suppose given a finite, flat, commutative group scheme $G$ over
$\Spec(F)$, and we let $G^* = \Spec(A)$. Let $\fA$ denote any
$O_F$-algebra in $A$ such that $\fA \otimes_{O_F} F = A$. Set
\begin{equation*}
\bH(\fA) := \{ \alpha \in \fA_{R^c}^{\times}\,|\, \alpha^{\omega} = g_{\omega}
\alpha \,\, \text{for all $\omega \in \Omega_K$, where $g_{\omega} \in
\Gamma$}\}.
\end{equation*}

We shall be interested in using the groups $\bH(\fA)$ and $\bH(A)$ to
impose Selmer-type conditions on elements of $H^1(F,G)$. The following
definition is motivated by Corollary \ref{C:Leon}(a). (Recall that the
isomorphism $\Up_F$ below was defined in Corollary \ref{C:Leon}.)

\begin{definition} \label{D:locsel} 
Suppose that $F$ is a local field. Then we define the subgroup
$H^1_{\fA}(F,G)$ of $H^1(F,G)$ 
by:
\begin{equation*}
H^1_{\fA}(F,G) = \left\{ x \in H^1(F,G) \mid \Up_F(x) \in
\frac{\bH(\fA) \cdot A^{\times}}{\Gamma \cdot A^{\times}} \subseteq
\frac{\bH(A)}{\Gamma \cdot A^{\times}}= H(A) \right\}.
\end{equation*}
Hence a $G$-torsor $\pi: X \to \Spec(F)$ lies in $H^1_{\fA}(F,G)$ if and
only if there exists a trivialisation $s_{\pi}:A \xrightarrow{\sim}
\L_{\pi}$ with $\br(s_{\pi}) \in \bH(\fA) \subseteq \bH(A)$. The
resolvend $\br(s_{\pi})$ of such a trivialisation is well-defined up
to multiplication by an element of $\Gamma \cdot \fA^{\times}$.
\qed
\end{definition}

\begin{definition} \label{D:glosel} 
If $F$ is a number field, then we define $H^1_{\fA}(F,G)$ by
\begin{equation*}
H^1_{\fA}(F,G) = \Ker \left[ H^1(F,G) \to \prod_{v \nmid \infty}
\frac{H^1(F_v,G)}{H^1_{\fA_v}(F_v,G)} \right].
\end{equation*}
\qed
\end{definition}

\begin{remark} \label{R:maxeta}
Suppose that $F$ is either a number field or a local field, and let
$\fM$ be the unique maximal $O_F$-order in $A$. Assume that $\pi \in
H^1(F,G)$ lies in $\Ker(\eta_F)$. Then it follows from the discussion
in Remark \ref{R:kereta} that there exists a resolvend $\br(s_\pi) \in
\bH(A)$ associated to $\pi$ such that $\br(s_\pi)^N = 1$. Hence
$\br(s_\pi) \in \bH(\fM)$, and so $\pi \in H^{1}_{\fM}(F,G)$. \qed
\end{remark}

Now suppose that $F$ is a number field, and let $\fM$ be the unique
maximal $O_F$-order in $A$. Let $J_f(A)$ denote the group of finite
ideles of $A$, i.e. $J_f(A)$ is the restricted direct product of the
groups $A_{v}^{\times}$ with respect to the subgroups
$\fM_{v}^{\times}$ for $v \nmid \infty$. We view $A^{\times}$ as being
a subgroup of $J_f(A)$ via the obvious diagonal embedding. Write
$\Cl(\fA)$ for the locally free classgroup of $\fA$. Thus, $\Cl(\fA)$
is the Grothendieck group of locally free $\fA$-modules of finite
rank, and it may be identified with the group $\Pic(\Sp(\fA))$. Then
it is a standard result from the theory of classgroups (see
e.g. \cite[\S52]{CR}) that there is a natural isomorphism
\begin{equation} \label{E:ali iso}
\Cl(\fA) \simeq \frac{J_f(A)}{\left(\prod_{v \nmid \infty}
\fA_{v}^{\times} \right) \cdot A^{\times}}.
\end{equation}

\begin{theorem} \label{T:classmap}
Let $F$ be a number field. 

(a)There is a natural homomorphism
\begin{equation*}
\phi_{\fA}: H^{1}_{\fA}(F,G) \to \Cl(\fA).
\end{equation*}

(b) The isomorphism
\begin{equation*}
\Up_F: H^1(F,G) \xrightarrow{\sim} H(A)
\end{equation*}
of Corollary \ref{C:Leon} induces an isomorphism
\begin{equation*}
\Up_{F,\fA}: \Ker(\phi_{\fA}) \xrightarrow{\sim}
H(\fA) \subseteq H(A).
\end{equation*}

(c) We have $\Ker(\eta_F) \subseteq \Ker(\phi_\fM)$.
\end{theorem}

\begin{proof}
(a) Suppose that $\pi:X \to \Spec(F)$ is a $G$-torsor with $[\pi] \in
H^1_{\fA}(F,G)$. Fix a trivialisation $s_{\pi}: A \xrightarrow{\sim}
\L_{\pi}$. Note that $\br(s_{\pi}) \in \bH(A)$ is well-defined up to
multiplication by an element of $\Gamma \cdot A^{\times}$.

For each finite place $v$ of $F$, write $\pi_v$ for the torsor $X
\otimes_F F_v \to \Spec(F_v)$. Since $[\pi_v] \in
H^{1}_{\fA_v}(F_v,G)$, we may choose a trivialisation $t_{\pi_v}: A_v
\xrightarrow{\sim} \L_{\pi_v}$ satisfying $\br(t_{\pi_v}) \in
\bH(\fA_v)$. Then $\br(t_{\pi_v})$ is well-defined up to
multiplication by an element of $\Gamma \cdot \fA^{\times}_{v}$, and
$\br(t_{\pi_v}) \br(s_{\pi}^{-1}) \in \Gamma \cdot
A_v^{\times}$. Furthermore, for all but finitely many places $v$, we
have that $\br(s_{\pi}) \in \bH(\fM_v)$, and so $\br(t_{\pi_v})
\br(s_{\pi}^{-1}) \in \Gamma \cdot \fM_v^{\times}$ for all such $v$.

It therefore follows that the element $(\br(t_{\pi_v})
\br(s_{\pi}^{-1}))_v$ lies in $\left(\prod_{v \nmid \infty} \Gamma
\right) \cdot J_f(A)$, and that its image in
\begin{equation*}
\frac{\left( \prod_{v \nmid \infty} \Gamma \right) \cdot
J_f(A)}{\left(\prod_{v \nmid \infty} \Gamma \right) \cdot \left
( \prod_{v \nmid \infty} \fA_v^{\times} \right) \cdot A^{\times}}
\simeq
\frac{J_f(A)}{\left(\prod_{v \nmid \infty}
\fA_{v}^{\times} \right) \cdot A^{\times}}
\simeq
\Cl(\fA)
\end{equation*}
is well-defined. We define 
\begin{equation*}
\phi_{\fA}(\pi) = [(\br(t_{\pi_v})
\br(s_{\pi}^{-1}))_v] \in \Cl(\fA).
\end{equation*}

We now show that $\phi_{\fA}$ is a homomorphism. Suppose that $\pi':X'
\to \Spec(F)$ is another $G$-torsor, and let
$s_{\pi'}$ and $t_{\pi'_v}$ ($v \nmid \infty$) be trivialisations
defined analogously to $s_{\pi}$ and $t_{\pi_v}$ above. Write $\pi'':=
\pi \cdot \pi'$. Then it follows from the functoriality of
Waterhouse's construction in \cite{W} that there is a natural
isomorphism $\L_{\pi''} \simeq \L_{\pi} \otimes \L_{\pi'}$. Thus, if
we set
\begin{equation*}
s_{\pi''}:= s_{\pi} \otimes s_{\pi'}: \L_{\pi} \otimes \L_{\pi'}
\xrightarrow{\sim} A, \qquad
t_{\pi_{v}''}:= t_{\pi_v} \otimes
t_{\pi_{v}'}: \L_{\pi_v} \otimes \L_{\pi_{v}'} \xrightarrow{\sim}
A_v,
\end{equation*} 
then $\br(s_{\pi''})= \br(s_{\pi}) \br(s_{\pi'})$, and
$\br(t_{\pi_{v}''}) = \br(t_{\pi_v}) \br(t_{\pi_{v}'})$. Hence it
follows that
\begin{align*}
\phi_{\fA}(\pi'') = [(\br(t_{\pi_{v}''}) \br(s_{\pi''})^{-1})_v] 
&= [(\br(t_{\pi_{v}})\br(s_{\pi}^{-1}))_v] \cdot
[(\br(t_{\pi'_{v}})\br(s_{\pi'}^{-1}))_v] \\
&= \phi_{\fA}(\pi) \cdot \phi_{\fA}(\pi'),
\end{align*}
as asserted.
\smallskip

(b) Suppose that $\pi:X \to \Spec(F)$ is a $G$-torsor satisfying
$\br(s_{\pi}) \in \bH(\fA)$ for some choice of trivialisation
$s_{\pi}: A \xrightarrow{\sim} \cL_{\pi}$ of $\cL_{\pi}$. Write
$s_{\pi,v}: A_v \xrightarrow{\sim} \cL_{\pi_v}$ for the trivialisation
of $\cL_{\pi_v}$ induced by $s_{\pi}$. Then $\br(s_{\pi,v}) \in
\bH(\fA_v)$ for all finite places $v$ of $F$. Hence $\pi \in
H^1_{\fA}(F,G)$, and we may take $t_{\pi_v} = s_{\pi,v}$ in the
definition of $\phi_{\fA}(\pi)$ given in part (a). This in turn
gives $\phi_{\fA}(\pi) = 0$.

Now suppose conversely that $\pi \in H^1_{\fA}(F,G)$ with
$\phi_{\fA}(\pi) = 0$. Then for suitable choices of $s_{\pi}$ and
$\br(t_{\pi_v})$, we have that
\begin{equation*}
(\br(t_{\pi_v}) \br(s_{\pi})^{-1})_v \in \prod_{v \nmid \infty}
\fA_v^{\times} \subseteq \prod_{v \nmid \infty} \bH(\fA_v).
\end{equation*}
This implies that $\br(s_{\pi}) \in \bH(\fA_v)$ for each place $v
\nmid \infty$, and so it follows that $\br(s_{\pi}) \in \bH(\fA)$.
\smallskip

(c) This follows directly from Remark \ref{R:maxeta} and part (b) above.
\end{proof}

\begin{example} \label{e:flat}
Suppose that $F$ is a number field. Let $\cG$ be a finite, flat,
commutative group scheme over $\Spec(O_F)$, with generic fibre
$G$. Let $\cG^* = \Spec(\fA)$ denote the Cartier dual of $\cG^*$; then
$G^* = \Spec(A)$ is the generic fibre of $\cG^*$. Corollary
\ref{C:Leon}(a) implies that $H^{1}_{\A_v}(F_v,G) = H^1(O_{F_v},G)$
for each finite place $v$ of $F$, and so it follows that
\begin{equation*}
H^1_{\fA}(F,G) = H^1(O_F, \cG).
\end{equation*}
Hence we obtain a description of the flat Selmer group of $\cG$ in terms of
resolvends. In this case, the map
\begin{equation*}
\phi_{\fA}: H^1(O_F,\cG) \to \Cl(\fA) \simeq \Pic(\cG^*)
\end{equation*}
is the same as the class invariant homomorphism \eqref{E:psi} for the
group $\cG$.

Also, we have
\begin{equation*}
H^1(O_F, \cG) = H^1_{\fA}(F,G) \subseteq H^1_{\fM}(F,G),
\end{equation*}
and $\Sp(\fM)$ is the normalisation of $\cG^*$. The restriction of the
homomorphism
\begin{equation*}
\phi_{\fM}: H^1_{\fM}(F,G) \to \Cl(\fM) \simeq \Pic(\Sp(\fM))
\end{equation*}
to $H^1(O_F,\cG)$ is the same as the class invariant homomorphism
\eqref{E:varphi}. \qed
\end{example}

\begin{remark} \label{R:unram}
Suppose that $F$ is a number field, and let $N$ denote the exponent of
$G$. If $v$ is a place of $F$ with $v \nmid N$, set
$$
H^1_f(F_v,G):= \Ker \left[ H^1(F_v,G) \to H^1(F^\nr_v,G) \right],
$$
where $F^\nr_v$ is the maximal unramified extension of $F_v$ in a
fixed algebraic closure of $F_v$. 

If $\pi \in H^1_f(F_v,G)$, and $\br(s_\pi)$ is any resolvend asociated
to $\pi$, then Proposition \ref{P:incres} implies that $\br(s_\pi) \in
A_{v,F^\nr_v}^{\times}$. Since $F^\nr_v/F_v$ is unramified, it follows
(via considering the Wedderburn decomposition \eqref{E:wedd} of $A_v$)
that there exists $\alpha \in A_v^\times$ such that $\alpha^{-1}
\br(s_\pi) = \br(\alpha^{-1} s_\pi) \in
\fM_{v,O_{F^\nr_v}}^{\times}$. This implies that $\pi \in
H^{1}_{\fM}(F,G)$, and so
$$
H^1_f(F_v,G) \subseteq H^{1}_{\fM}(F_v,G).
$$

Suppose further that $G$ is unramified at $v$. Then $\cG:= \Sp(\fM_v)$
is a finite, flat, commutative $O_{F_v}$-group scheme, and it is a
standard result that $H^1_f(F_v,G) = H^1(O_{F_v},\cG)$. We therefore
deduce that in this case, we have $H^1_f(F_v,G) = H^{1}_{\fM}(F_v,G)$.
\qed
\end{remark}

\begin{remark} \label{R:twisted}
It is not difficult to define refinements of the homomorphism
$\phi_\fA$ taking values in relative algebraic $K$-groups as in
\cite{AB}, or in Arakelov Picard groups as in \cite{AP}. However, for
the sake of brevity, we shall not go into this here. \qed
\end{remark}

Now suppose that $F$ is a number field, and let $L/F$ be a finite
extension. It is not hard to check that the homomorphism $\cN_{L/F}$
of \eqref{E:norm} induces a homomorphism
$$
\cN_{L/K}: \Cl(\fA_{O_L}) \to \Cl(\fA).
$$

\begin{proposition} \label{P:classnorm}
If $F$ is a number field, and $L/F$ is a finite extension, then the
following diagram is commutative:
\begin{equation} \label{E:classnorm}
\begin{CD}
H^{1}_{\fA_{O_L}}(F,G) @>{\phi_{\fA_{O_L}}}>>  \Cl(\fA_{O_L})  \\
@V{\Cores_{L/F}}VV                   @V{\cN_{L/F}}VV \\
H^{1}_{\fA}(F,G) @>{\phi_{\fA}}>>   \Cl(\fA).
\end{CD}
\end{equation}
\end{proposition}

\begin{proof} Let $\pi:X \to \Spec(L)$ be a $G$-torsor with $[\pi] \in
H^{1}_{\fA_{O_L}}(L,G)$, and let $s_\pi: A_L \xrightarrow{\sim}
\cL_\pi$ be any trivialisation of $\cL_\pi$. For each finite place $v$
of $L$, let $t_{\pi_v}: A_{L_v} \xrightarrow{\sim} \cL_{\pi_v}$ be a
trivialisation of $\cL_{\pi_v}$ satisfying $\br(t_{\pi_v}) \in
\bH(\fA_{O_{L_v}})$, Then
$$
\phi_{\fA_{O_L}}(\pi) = [(\br(t_{\pi_v})\br(s_\pi)^{-1})_v] \in
\Cl(\fA).
$$
The result now follows via a similar argument to that used in the
proof of Theorem \ref{T:traceres}.
\end{proof}

For the rest of this paper, we shall mainly be concerned with the
special cases in which $\fA=\fM$ or $\fA= \fM \otimes_{O_F} O_F[1/p]:=
\fM^\brp$. We identify $\Pic(\Sp(\fM))$ and $\Pic(\Sp(\fM^\brp))$ with
the locally free classgroups $\Cl(\fM)$ and $\Cl(\fM^\brp)$ of $\fM$
and $\fM^\brp$ respectively.  We set
\begin{equation*}
H^1_u(F,G):=H^1_{\fM}(F,G),\qquad H^{1}_{u,\brp}(F,G):=
H^{1}_{\fM^\brp}(F,G),
\end{equation*}
and we write
\begin{equation*}
\phi: H^1_u(F,G) \to \Cl(\fM),\quad
\phi^\brp:H^{1}_{u,\brp}(F,G) \to \Cl(\fM^\brp)
\end{equation*}
for the homomorphisms given by Theorem \ref{T:classmap}. 

\begin{proposition} \label{P:kerphi iso}
Let $F$ be a number field, and suppose that $G^*$ is constant over
$\Spec(F)$.

(i) If $v$ is any finite place of $F$, then the isomorphism
$\eta_{F_v}$ of Proposition \ref{P:etaprim} induces an isomorphism
\begin{equation*}
H^1_u(F_v,G) \xrightarrow{\sim} \Hom(\Gamma^*,
O_{F_v}^{\times}/(O_{F_v}^{\times})^N).
\end{equation*}

(ii) The isomorphism $\eta_F$ induces isomorphisms
\begin{align*}
&\Ker(\phi) \xrightarrow{\sim} \Hom(\Gamma^*,
O_{F}^{\times}/(O_{F}^{\times})^N),\\
&\Ker(\phi^\brp) \xrightarrow{\sim} \Hom(\Gamma^*,
O_{F}[1/p]^{\times}/(O_{F}[1/p]^{\times})^N)
\end{align*}

\end{proposition}

\begin{proof} Since $G^*$ is constant over $\Spec(F)$, we have 
\begin{equation*}
A \simeq \Map(\Gamma^*,F),\qquad 
\fM \simeq \Map(\Gamma^*,O_F), \qquad
\fM^\brp \simeq \Map(\Gamma^*,O_F[1/p]).
\end{equation*}

Proposition \ref{P:etaprim} implies that we have isomorphisms
\begin{align} 
&H^1(F,G) \simeq \Hom(\Gamma^*, F^{\times}/(F^{\times})^N),
\label{E:i} \\
&H^1(F_v,G) \simeq \Hom(\Gamma^*,
F_v^{\times}/(F_v^{\times})^N). \label{E:ii}
\end{align}

Hence (i) follows from \eqref{E:ii} and the definition of
$H^1(F_v,G)$, while (ii) follows from \eqref{E:i} together with Theorem
\ref{T:classmap}(b).
\end{proof}



\section{$p$-adic representations}


In this section, we shall apply our previous work to the situation
described in the introduction. We first recall the relevant notation.

Let $F$ be a number field and $V$ be a $d$-dimensional $\Q_p$-vector
space. Suppose that $\rho: \Omega_F \to \GL(V)$ is a continuous
representation which is ramified at only finitely many primes of
$F$. We set $V^*:= \Hom_{\Q_p}(V,\Q_p(1))$, and we write $\rho^*:
\Omega_F \to \GL(V^*)$ for the corresponding representation of
$\Omega_F$. Let $T \subseteq V$ be any $\Omega_F$-stable lattice, and
write $T^*:= \Hom_{\Z_p}(T, \Z_p(1))$. For each positive integer $n$,
we define finite group schemes $G_n$ and $G_n^*$ over $\Spec(F)$ by
\begin{equation*}
G_n(F^c):= \Gamma_n= p^{-n}T/T; \qquad G_n^*(F^c):= \Gamma_n^*=
p^{-n}T^*/T^*.
\end{equation*} 
Then $G_n^*$ is the Cartier dual of $G_n$ with $G_n^* =
\Spec(A_n)$ for the Hopf algebra $A_n = (F^c \Gamma_n)^{\Omega_F}$
over $F$. 

Recall that $q_n: A_n \to A_{n-1}$ is the homomorphism induced by the
dual $p^D_n$ of the multiplication-by-$p$ map $p_n:G_n \to
G_{n-1}$. Suppose that, for each $n$, we are given an $O_F$-algebra
$\fA_n = \subseteq A_n$ such that $\fA_n \otimes_{O_F} F =
A_n$, and $q_n(\fA_n) = \fA_{n-1}$. Then it is easy to check that
$q_n$ induces homomorphisms
$$
\bH(A_n) \to \bH(A_{n-1}),\quad \text{and}\quad \bH(\fA_{n,v}) \to
\bH(\fA_{n-1,v})
$$
for each finite place $v$ of $F$. This implies that the natural maps
$$
H^1(F_v,G_n) \to H^1(F_v,G_{n-1}),\quad H^1(F,G_n) \to H^1(F,G_{n-1})
$$
induce homomorphisms
$$
H^{1}_{\fA_{n,v}}(F_v,G_n) \to H^{1}_{\fA_{n-1,v}}(F_v,G_{n-1}),\quad
H^{1}_{\fA_{n}}(F,G_n) \to H^{1}_{\fA_{n-1}}(F,G_{n-1})
$$
via restriction.

Set $\fA(T):= \varprojlim \fA_n$ and $\fA_v(T):= \varprojlim
\fA_{n,v}$ (where the inverse limits are taken with respect to the
maps $q_n$), and let
\begin{align*}
\bH(\fA_v(T)) &:= \{ \alpha \in {\fA_{v}(T)}_{O_{F_v^c}}^{\times}\,|\,
\alpha^{\omega} = t_{\omega} 
\alpha \,\, \text{for all $\omega \in \Omega_{F^c_v}$, where $t_{\omega} \in
T$}\}, \\
H(\fA_v(T))&:= \frac{\bH(\fA_v(T))}{T \cdot {\fA_{v}(T)}^{\times}}. 
\end{align*}
Define $\bH(\fA(T))$ and $H(\fA(T))$ in a similar way.
Write
$$
H^{1}_{\fA_v(T)}(F_v,T):= \varprojlim H^{1}_{\fA_{n,v}}(F_v,G_n),\quad
H^{1}_{\fA(T)}(F,T):= \varprojlim H^{1}_{\fA_n}(F,G_n).
$$

\begin{proposition} \label{P:invleon}
For each finite place $v$ of $F$, we have
$$
H^{1}_{\fA_v(T)}(F_v,T) \simeq \frac{\bH(\fA_v(T))}{T \cdot
\fA_v^\times(T)}.
$$
\end{proposition}

\begin{proof}
It follows from the definition of $H^{1}_{\fA_{n,v}}(F_v,G_n)$ (see
Definition \ref{D:locsel}) that, for each $n$, there is an exact
sequence
$$
1 \to G_n \cdot \fA^{\times}_{n,v} \to \bH(\fA_{n,v}) \to
H^{1}_{\fA_{n,v}}(F_v,G_n) \to 0.
$$
Passing to inverse limits, and using the fact that the inverse system
$\{G_n \cdot \fA_{n,v}^{\times} \}_{n}$ satisfies the Mittag-Leffler
condition yields
$$
H^{1}_{\fA_{v}(T)}(F_v,T) \simeq \frac{\varprojlim \bH(\fA_{n,v})}{T
\cdot \fA_{v}(T)^{\times}}.
$$
It follows easily from the definitions that
$$
\bH(\fA_v(T)) = \varprojlim \bH(\fA_{n,v}),
$$
and this implies the result.
\end{proof}

It is easy to check that $p_n^D$ induces pullback homomorphisms
$$
(p_n^D)^*:\Cl(\fA_n) \to \Cl(\fA_{n-1}).
$$
Let
$$
\phi_{\fA_n}: H^{1}_{\fA_n}(F,G_n) \to \Cl(\fA_n)
$$
denote the natural homomorphism afforded by Theorem \ref{T:classmap}.

\begin{theorem} \label{T:limitmaps}
The following diagram is commutative:
\begin{equation} \label{E:compatible}
\begin{CD}
H^{1}_{\fA_n}(F,G_n) @>{\phi_{\fA_n}}>> \Cl(\fA_n) \\
@V{p_n}VV                     @VV{(p_{n}^{D})^{*}}V \\
H^{1}_{\fA_n}(F,G_{n-1}) @>{\phi_{\fA_{n-1}}}>>
\Cl(\fA_{n-1}), 
\end{CD}
\end{equation}
\end{theorem}

\begin{proof} The proof of this is similar to that of Theorem
\ref{T:pdiv}. Let $\pi_n:X_n \to \Spec(F)$ be a $G_n$-torsor with $[\pi_n] \in
H^{1}_{\fA_n}(F,G_n)$, and write $\pi_{n-1}:= p_n(\pi_n)$. Let $s_{\pi_n}: A_n
\xrightarrow{\sim} \L_{\pi_n}$ be any trivialisation of $\L_{\pi_n}$,
and for each finite place $v$ of $F$, let $t_{\pi_{n},v}: A_{{n},v}
\xrightarrow{\sim} \L_{{\pi_n},v}$ be a trivialisation of
$\L_{{\pi_n},v}$ satisfying $\br(t_{\pi_{n},v}) \in \bH(\fA_{n,v})$.

Then, via functoriality, we have that $(p_n^D)^*\L_{\pi_n} \simeq
\L_{\pi_{n-1}}$. The pullbacks $(p_n^D)^* s_{\pi_n}$ and
$(p_n^D)^*t_{\pi_{n},v}$ of $s_{\pi_n}$ and $t_{\pi_{n},v}$ along $p_n^D$
give trivialisations of $\L_{\pi_{n-1}}$ and $\L_{{\pi_{n-1}},v}$
respectively, and we have that
\begin{equation*}
\br((p_n^D)^*t_{\pi_{n},v}) = q_{n,v}(t_{\pi_{n},v}) \in \bH(\fA_{n-1,v}).
\end{equation*}
The result now follows from the definitions of $\phi_{\fA_n}$
and $\phi_{\fA_{n-1}}$.
\end{proof}

Set
$$
\Cl(\fA(T)) := \varprojlim \Cl(\fA_n).
$$
Then passing to inverse limits over the diagrams \eqref{E:compatible}
yields a homomorphism
\begin{equation} \label{E:picmap}
\Phi_{\fA(T)}: H^{1}_{\fA(T)}(F,T) \to \Cl(\fA(T)).
\end{equation}

\begin{proposition} \label{P:infclassnorm}
Suppose that $L/F$ is a finite extension. Then the map $\cN_{L/K}$
(see \eqref{E:norm} induces a homomorphism
$$
\cN_{L/K}: \Cl(\fA(T)_{O_L}) \to \Cl(\fA(T)),
$$
and the following diagram is commutative:
\begin{equation} \label{E:infclassnorm}
\begin{CD}
H^{1}_{\fA(T)_{O_L}}(F,G) @>{\phi_{\fA(T)_{O_L}}}>>  \Cl(\fA(T)_{O_L})  \\
@V{\Cores_{L/F}}VV                   @V{\cN_{L/F}}VV \\
H^{1}_{\fA(T)}(F,G) @>{\phi_{\fA(T)}}>>   \Cl(\fA(T)).
\end{CD}
\end{equation}
\end{proposition}

\begin{proof} This follows from Theorem \ref{P:classnorm}.
\end{proof}

\begin{proposition}
There is an isomorphism
$$
\Upsilon_{F,\fA(T)}: \Ker(\Phi_{\fA(T)}) \xrightarrow{\sim} H(\fA(T))
\subseteq H(A(T)).
$$
\end{proposition}

\begin{proof} This follows easily from Theorem \ref{T:classmap}(b).
\end{proof}

Write 
\begin{equation*}
H^1_u(F,T):= \varprojlim H^1_u(F,G_n),\qquad 
H^{1}_{u,\brp}(F,T):= \varprojlim H^{1}_{u,\brp}(F,G_n).
\end{equation*} 
Then \eqref{E:picmap} yields a homomorphism
\begin{equation}
\Phi_F:= \Phi_{\fM(T)}: H^1_u(F,T) \to \varprojlim \Cl(\fM(T)).
\end{equation}
From Remark \ref{R:unram}, we see that
$$
H^{1}_{f,\brp}(F,T) \subseteq H^{1}_{u,\brp}(F,T).
$$
Hence, restricting $\Phi_{\fM^\brp(T)}$ to $H^{1}_{f,\brp}(F,T)$ yields a
homomorphism
$$
\Phi^\brp_F: H^{1}_{f,\brp}(F,T) \to \varprojlim \Cl(\fM^\brp(T)).
$$

\begin{remark} \label{R:bloch-kato}
If $v \mid p$, then $H^1_u(F_v,T)$ is not in general equal to the
group $H^1_f(F_v,T)$ introduced by Bloch and Kato in \cite{BK}. For
example, it follows from the definitions that $H^1_u(F_v,T)$ is always
infinite. On the other hand, if $v \mid p$, then $H^1_f(F_v,\Z_p(r))$
is finite for all $r <0$ (see e.g. \cite[Example 3.9]{BK}). 

Let $S$ be any finite set of places of $F$ containing all places lying
above $p$, as well as all places at which $T$ is ramified, and let
$F^S/F$ denote the maximal extension of $F$ which is unramified
outside $S$. Then it follows from the definitions that $H^1_u(F,T)
\subseteq H^1(F^S/F,T)$, and so we deduce that $H^1_u(F,T)$ is always
a finitely generated $\bZ_p$-module. \qed

\end{remark}

\begin{remark} \label{R:conjecture}
Suppose that $\cA$ is an abelian scheme over $\Spec(O_F)$, and let $T$
denote its $p$-adic Tate module. For each positive integer $n$, let
$\cG_n$ denote the $O_F$-group scheme of $p^n$-torsion on $\cA$, and
write $\cG_n^*$ for its Cartier dual. Then taking inverse limits of
the homomorphisms
\begin{equation*}
\psi_n: H^1(\Spec(O_F),\cG_n) \to \Pic(\cG_n^*)
\end{equation*}
yield a homomorphism
\begin{equation*}
\Psi_F: H^1_f(F,T) \to \varprojlim \Pic(\cG_n^*)
\end{equation*}
(see \cite{A2}, \cite{AT}, \cite{AP}). It seems reasonable to
conjecture that $\Psi_F$ is injective modulo torsion. In \cite{AT},
this conjecture is shown to be true (subject to certain technical
hypotheses) when $\cA_{/F}$ is an elliptic curve and $p$ is a prime of
ordinary reduction. \qed
\end{remark}


\section{Proof of Theorem \ref{T:intro 2}} \label{S:max}


In this section we give the proof of Theorem \ref{T:intro 2}.
\smallskip

For each integer $n$ the action of $\Omega_F$ on $\Gamma_n^*$ yields a
representation
\begin{equation*}
\rho_n^*: \Omega_F \to \Aut(\Gamma_n^*).
\end{equation*}
Write $F_n^*$ for the fixed field of $\rho_n^*$; then $F_{\infty}^* =
\cup_n F_n^*$, where $F_{\infty}^*$ is the extension of $F$ cut out by
$\rho^*$. The group scheme $G_n^*$ is constant over $\Spec(F_n^*)$,
and we write
\begin{equation*}
\eta_{n,F_n^*}: \Ker(\phi_{n,F_n^*}) \xrightarrow{\sim} \Hom(\Gamma_n^*,
O_{F_n^*}^{\times}/ (O_{F_n^*}^{\times})^{p^n})
\end{equation*}
for the isomorphism afforded by Proposition \ref{P:kerphi iso}(ii).

Consider the map
\begin{equation*}
d_n: H^1(F^*_n,\Gamma_n) \xrightarrow{p_n} H^1(F^*_n, \Gamma_{n-1})
\xrightarrow{\Cores_{F^*_n/F^{*}_{n-1}}} H^1(F^{*}_{n-1}, \Gamma_{n-1}).
\end{equation*}

\begin{lemma} \label{L:rubin} Passing to the inverse limit of the maps
\begin{equation*}
d_n: H^1(F^*_n, \Gamma_n) \to H^1(F^{*}_{n-1}, \Gamma_{n-1})
\end{equation*}
induces isomorphisms
\begin{align}
\varprojlim H^1(F^*_n, \Gamma_n) &\xrightarrow{\sim} \varprojlim
H^1(F^*_n,T), \label{E:rubin}  \\
\varprojlim H^1(F^{*}_{n,v}, \Gamma_n) &\xrightarrow{\sim} \varprojlim
H^1(F^{*}_{n,v},T) \label{E:rubinlocal}
\end{align}
\end{lemma}

\begin{proof} See e.g. \cite[Lemma B.3.1]{R}.
\end{proof}

Let $h_n$ denote the composition
\begin{align*}
h_n: \Hom(\Gamma_n^*, F_{n}^{*\times}/(F_{n}^{*\times})^{p^n})
&\xrightarrow{r_n}
\Hom(\Gamma_{n-1}^*, F_{n}^{*\times}/(F_{n}^{*\times})^{p^{n-1}}) \\
&\xrightarrow{N_{F^{*}_{n}/F^{*}_{n-1}}}
\Hom(\Gamma_{n-1}^*, F_{n-1}^{*\times}/(F_{n-1}^{*\times})^{p^{n-1}}),
\end{align*}
where $r_n$ is defined in Corollary \ref{C:pdiv compatible}, and
$N_{F^{*}_{n}/F^{*}_{n-1}}$ is induced by the norm map from $F^*_n$ to
$F^{*}_{n-1}$.

\begin{lemma} \label{L:bpr}
Passing to the inverse limit of the maps
\begin{equation*}
h_n: \Hom(\Gamma_n^*, F_{n}^{*\times}/(F_{n}^{*\times})^{p^n}) \to
\Hom(\Gamma_{n-1}^*, F_{n-1}^{*\times}/(F_{n-1}^{*\times})^{p^{n-1}})
\end{equation*}
induces isomorphisms
\begin{equation*}
\varprojlim \Hom(\Gamma_n^*, F_{n}^{*\times}/(F_{n}^{*\times})^{p^n})
\xrightarrow{\sim}
\Hom(T^*, \varprojlim \Check{F}_{n}^{*\times}),
\end{equation*}
\begin{equation*}
\varprojlim \Hom(\Gamma_n^*,
O_{F^{*}_{n}}^{\times}/(O_{F^{*}_{n}}^{\times})^{p^n})
\xrightarrow{\sim}
\Hom(T^*, \Check{O}_{F^*_n}^{\times}).
\end{equation*}
\end{lemma}

\begin{proof} This is proved by applying Lemma \ref{L:rubin} to the
$p$-divisible group schemes $(\Z /p^n \Z)_{n \geq 1}$ and
$(\boldsymbol{\mu}_{p^n})_{n \geq 1}$. We have
\begin{align*}
\varprojlim[\Hom(\Gamma_n^*, F_{n}^{*\times}/(F_{n}^{*\times})^{p^n})]
&=
\varprojlim[\Hom(\Gamma_n^*,\Z/p^n\Z) \otimes_{\Z_p}
(F_{n}^{*\times}/(F_{n}^{*\times})^{p^n})] \\
&=
[\varprojlim \Hom(\Gamma_n^*,\Z/p^n\Z)] \otimes_{\Z_p} [\varprojlim
(F_{n}^{*\times}/(F_{n}^{*\times})^{p^n})] \\
&=
\Hom(T^*, \Z_p) \otimes_{\Z_p} \varprojlim H^1(F^*_n, \mu_{p^n}) \\
&=
\Hom(T^*, \Z_p) \otimes_{\Z_p} \varprojlim H^1(F^*_n,\Z_p(1)) \\
&=
\Hom(T^*, \Z_p) \otimes_{\Z_p} \varprojlim \Check{F}_{n}^{*\times} \\
&=
\Hom(T^*, \varprojlim \Check{F}_{n}^{*\times}).
\end{align*}
The second isomorphism may be established in a similar manner.
\end{proof}

\begin{lemma}
The following diagram is commutative:
\begin{equation} \label{E:key compatible}
\begin{CD}
H^1(F^*_n,\Gamma_n) @>{\eta_{n,F^*_n}}>>
\Hom(\Gamma_n^*,F_{n}^{*\times}/(F_{n}^{*\times})^{p^n}) \\
@V{d_n}VV                     @VV{h_n}V              \\
H^1(F^{*}_{n-1},\Gamma_{n-1}) @>{\eta_{n,F^{*}_{n-1}}}>>
\Hom(\Gamma_{n-1}^{*},F_{n-1}^{*\times}/(F_{n-1}^{*\times})^{p^{n-1}}).
\end{CD}
\end{equation}
\end{lemma}

\begin{proof} Note that $G_n^*$ is constant over $\Spec(F^*_n)$. The
result now follows from Theorem \ref{T:traceres} and Corollary
\ref{C:pdiv compatible}.
\end{proof}

Passing to inverse limits over the diagrams \eqref{E:key compatible},
and applying Lemmas \ref{L:rubin} and \ref{L:bpr} yields a natural
isomorphism
\begin{equation} \label{E:great iso}
\beta: \varprojlim H^1(F^*_n,T) \xrightarrow{\sim} \Hom(T^*, \varprojlim
\Check{F}_{n}^{*\times}).
\end{equation}
For each finite place $v$ of $F$, similar arguments to those given
above show that there is also a local isomorphism
\begin{equation} \label{E:great local iso}
\beta_v:\varprojlim H^1(F^{*}_{n,v},T) \xrightarrow{\sim} \Hom(T^*, \varprojlim
\Check{F}_{n,v}^{*\times}).
\end{equation}

\begin{proposition} \label{P:beta}
(i) For each finite place $v$ of $F$, the map $\beta_v$ induces an
isomorphism (which we denote by the same symbol)
\begin{equation} \label{E:beta local}
\beta_v: \varprojlim H^1_u(F^{*}_{n,v},T) \xrightarrow{\sim} \Hom(T^*,
\varprojlim \Check{O}_{F^{*}_{n,v}}^{\times}).
\end{equation}

(ii) The map $\beta$ induces an isomorphism
\begin{equation} \label{E:beta global}
\beta: \varprojlim H^1_u(F^{*}_{n},T) \xrightarrow{\sim} \Hom(T^*,
\varprojlim \Check{O}_{F^*_n}^{\times}).
\end{equation}
\end{proposition}

\begin{proof}
(i) It is easy to check that \eqref{E:rubinlocal} induces an isomorphism
\begin{equation*}
\varprojlim H_u^1(F^{*}_{n,v}, \Gamma_n) \xrightarrow{\sim} \varprojlim
H_u^1(F^{*}_{n,v},T).
\end{equation*}
The result now follows from the isomorphism
\begin{equation*}
H_u^1(F^{*}_{n,v},G_n) \xrightarrow{\sim}
\Hom(\Gamma_n^*, O_{F^{*}_{n,v}}^{\times}/(O_{F^{*}_{n,v}}^{\times})^{p^n})
\end{equation*}
afforded by Proposition \ref{P:kerphi iso}(i).

(ii) Suppose that $f \in \Hom(T^*, \varprojlim
\Check{F}_n^{*\times})$. For each finite place $v$ of $F$, write $f_v$
for the image of $f$ in $\Hom(T^*, \varprojlim
\Check{F}_{n,v}^{*\times})$. Then $f_v \in \Hom(T^*,
\varprojlim{O}_{F^*_n,v}^{\times})$ for all $v$ if and only if $f \in
\Hom(T^*, \varprojlim \Check{O}_{F^*_n}^{\times})$. The result now
follows from (i) above.
\end{proof}

\begin{corollary} \label{C:univ ker}
There is a natural isomorphism
\begin{equation*}
\varprojlim H_u^1(F^*_n,T) \simeq \varprojlim \Ker(\Phi_{F^*_n}).
\end{equation*}
\end{corollary}

\begin{proof} This follows directly from \eqref{E:beta local}, Lemma
\ref{L:bpr}, and Proposition \ref{P:kerphi iso}(ii).
\end{proof}

Corollary \ref{C:univ ker} implies that we have
$$
\cap_n \Cores_{F^{*}_{n}/F}(H^1_u(F^*_n,T)) \subseteq \Ker(\Phi_F).
$$
Hence, if $x \in \fC_F(T)$, then for some integer $M>0$, we have
$\Phi_F(Mx) = M \Phi_F(x) =0$, and so $\Phi_F(x)$ is of
finite order. This proves Theorem \ref{T:intro 2}.


\section{Proof of Theorem \ref{T:intro 1}} \label{S:pmax}


This section is devoted to the proof of Theorem \ref{T:intro 1}.
\smallskip

Recall that $C_{\infty} = \cup_n C_n$ denotes the cyclotomic
$\bZ_p$-extension of $F$. We begin by showing that if $x \in
\fG_F(T)$, then $\Phi_F^\brp(x)$ is of finite order.

\begin{lemma} \label{L:maxtriv}
Suppose that $m \geq 1$. Then
$$
\varprojlim \Cl(\fM^{\brp}_{m,O_{C_n}}) = 0,
$$
where the inverse limit is taken with respect to the maps
$\cN_{C_{n}/C_{n-1}}$.
\end{lemma}

\begin{proof}
In order to ease notation somewhat, we write $\fM^\brp$ for
$\fM^\brp_m$. Using the notation of Remark \ref{R:wedd}, we have
$$
\fM^{\brp}_{C_n} \simeq \prod_{\gamma^* \in \Gamma_m \backslash
\Omega_{C_n}} O_{C_n[\gamma^*]} [1/p],
$$
and so
$$
\Cl(\fM^{\brp}_{C_n}) \simeq \prod_{\gamma^* \in \Gamma_m \backslash
\Omega_{C_n}} \Cl(O_{C_n[\gamma^*]} [1/p]).
$$
Choose $d \geq 1$ sufficiently large that
$$
\Gamma_m^* \backslash \Omega_{C_d} = \Gamma^*_m \backslash
\Omega_{C_i}
$$ for all $i \geq d$. Then, for each $\gamma^* \in \Gamma^*_m
\backslash \Omega_{C_d}$
$$
C_\infty[\gamma^*]:= \cup_{i \geq d} C_i[\gamma^*]
$$
is the cyclotomic $\bZ_p$-extension of $C_d[\gamma^*]$, and the map
$$
\cN_{C_{i}/C_{d}}: \Cl(\fM^{\brp}_{O_{C_i}}) \to
\Cl(\fM^{\brp}_{C_d})
$$
is induced by the usual norm maps
$$
\Cl(O_{C_i[\gamma^*]}[1/p]) \to \Cl(O_{C_d[\gamma^*]}[1/p])
$$
for each $\gamma^* \in \Gamma^*_m \backslash \Omega_{C_i} =
\Gamma^*_m \backslash \Omega_{C_d}$.

However, if $L$ is any number field, and $L_\infty = \cup_j L_j$ is
the cyclotomic $\bZ_p$-extension of $L$, then it is easy to show that
$\varprojlim \Cl(O_{L_j}[1/p]) = 0$. Hence we have
$$
\varprojlim \Cl(O_{C_i[\gamma^*]}[1/p]) =0  \quad (i \geq d),
$$
and so it follows that
$$
\varprojlim \Cl(\fM^{\brp}_{O_{C_n}}) = 0,
$$
as claimed.
\end{proof}

\begin{corollary} \label{C:bigmaxtriv}
We have $\varprojlim \Cl(\fM^{\brp}_{O_{C_n}}(T)) = 0$, where the
inverse limit is taken with respect to the maps
$\cN_{C_{n}/C_{n-1}}$. \qed
\end{corollary}

\begin{proposition} \label{P:normker}
Suppose that $x \in \fG_F(T)$. Then $\Phi^\brp_F(x)$ is of finite
order.
\end{proposition}

\begin{proof} Taking inverse limits over \eqref{E:infclassnorm} yields
a homomorphism
$$
\varprojlim \Phi^{\brp}_{C_n}: \varprojlim
H^{1}_{\fM(T)_{O_{C_n}}}(C_n,T) \to \varprojlim
\Cl(\fM^{\brp}_{O_{C_n}}(T)),
$$
and Corollary \ref{C:bigmaxtriv} implies that this is the zero
map. Hence 
$$
\cap_n \Cores_{C_n/F}(H^{1}_{\fM(T)_{O_{C_n}}}(C_n,T)) \subseteq
\Ker(\Phi^\brp_F).
$$
Thus, if $x \in \fG_F(T)$, then $M \Phi_F^\brp(x)
=\Phi_F^\brp(Mx)= 0$ for some integer $M>0$. This establishes the
result.
\end{proof}

We now recall the definition of the pairing
$$
B_F: H^{1}_{f,\brp}(F,T) \times \Ker(\Loc_{F,T^*}) \to \bQ_p
$$
given in \cite[Section 3.1.4]{PR1}.

Fix $x \in H^{1}_{f,\brp}(F,T)$ and $y \in \Ker(\Loc_{F,T^*})$. Then
viewing $y$ as an element of $H^1(F,T^*) \simeq 
\Ext^{1}_{\Omega_F}(\bZ_p,T^*)$ yields an extension
\begin{equation} \label{E:ext1}
0 \to T^* \to T'_y \to \bZ_p \to 0.
\end{equation}
Taking $\bZ_p(1)$-duals of \eqref{E:ext1} yields an exact sequence
\begin{equation} \label{E:ext2}
1 \to \bZ_p(1) \to T_y \to T \to 0.
\end{equation}

We may consider the global and local Galois cohomology of
\eqref{E:ext2} for each finite place $v$ of $F$:
\begin{equation} \label{E:cohom}
\begin{CD}
H^1(F,\bZ_p(1)) @>{i}>> H^1(F,T_y) @>{j}>> H^1(F,T) @>>>
H^2(F,\bZ_p(1)) \\
@VVV  @VVV  @VVV  @VVV \\
H^1(F_v,\bZ_p(1)) @>{i_v}>> H^1(F_v,T_y) @>{j_v}>> H^1(F_v,T) @>>>
H^2(F_v,\bZ_p(1)).
\end{CD}
\end{equation}
It may be shown via Tate local duality that
$$
H^1_f(F_v,T) \subseteq j_v(H^1_f(F_v,T_y))
$$
for all places $v \nmid p$.

At places $v \mid p$ the extension \eqref{E:ext2} splits locally at
$v$ because $y \in \Ker(\Loc_{F,T^*})$, and so we have a corresponding
splitting 
\begin{equation} \label{E:csplit}
H^1(F_v,T_y) = H^1(F_v,\bZ_p(1)) \oplus H^1(F_v,T)
\end{equation}
on the level of cohomology groups. Hence we have 
$$
H^1_f(F_v,T) = j_v(H^1_f(F_v,T_y))
$$
in this case, and in fact every element $z \in H^1(F_v,T)$ has
a canonical lifting to an element of $H^1(F_v,T_y)$ given by $z
\mapsto (0,z)$.

Global classfield theory implies that the natural map
$$
H^2(F,\bZ_p(1)) \to \bigoplus_v H^2(F_v,\bZ_p(1))
$$
is injective, and so we deduce from \eqref{E:cohom} that
$$
H^{1}_{f,\brp}(F,T) \subseteq j(H^{1}_{f,\brp}(F,T_y)).
$$

Choose a global lifting $\ti{x} \in H^{1}_{f,\brp}(F,T_y)$ of $x \in
H^{1}_{f,\brp}(F,T)$. For each place $v$ with $ v \nmid p$, choose any
local lifting $\lambda_v \in H^1_f(F_v,T_y)$ of $x_v \in
H^1_f(F_v,T)$. For places $v$ with $v \mid p$, define $\lambda_v \in
H^1(F_v,T_y)$ to be the canonical lifting of $x_v$ afforded by the
splitting \eqref{E:csplit}. Then for each place $v$ of $F$, we have
$\ti{x}_v - \lambda_v \in i_v(H^1(F_f, \bZ_p(1)))$. If $v \mid p$,
then $i_v$ is injective, and so we may in fact identify
$i_v(H^1(F_v,\bZ_p(1)))$ with $H^1(F_v,\bZ_p(1))$.

Let 
$$
l_\chi: \bigoplus_v H^1(F_v, \bZ_p(1)) \to \bQ_p
$$
denote the composition
$$
\bigoplus_v H^1(F_v,\bZ_p(1)) \simeq \oplus_v \check{F}_v^\times
\xrightarrow{L_\chi} \bQ_p,
$$
where $L_\chi$ is defined by
$$
L_\chi((u_v)_v) = \sum_{v \mid p} \log_p N_{F_{v}/\bQ_{p}}(u_v) -
\sum_{v \nmid p} (\log_p q_v) \ord_v(u_v).
$$
(Here $q_v$ denotes the cardinality of the residue field of $F_v$, and
we choose Iwasawa's branch of the $p$-adic logarithm, so that
$\log_p(p)=0$.)

It may be shown that $l_\chi$ induces a well-defined map on $\oplus_v
i_v(H^1(F_v,\bZ_p(1)))$. We define
$$
B_F(x,y) = l_\chi \left( \ti{x}_v - \sum_v \lambda_v \right) \in \bQ_p.
$$
It is shown in \cite[Section 3.1.4]{PR1} that $B_F$ induces a pairing
$$
\langle \langle\,,\,\rangle \rangle:
\frac{H^1_{f,\{p\}}(F,T)}{\fG_F(T)} \times \Ker(\Loc_{F,T^*}) \to \Q_p,
$$
and it is conjectured that this pairing is always non-degenerate. We
shall relate this pairing to the homomorphism $\Phi_{\fM^\brp(T)}$ by
interpreting the pairing $B_F$ in terms of resolvends.

In order to ease notation in what follows, we shall write
$\Phi^{\brp}_{\bZ_p(1)}$, $\Phi^{\brp}_{T_y}$ and $\Phi^\brp_T$ for
$\Phi_{\fM^\brp(\bZ_p(1))}$, $\Phi_{\fM^\brp(T_y)}$, and
$\Phi_{\fM^\brp(T)}$ respectively.

\begin{proposition} \label{P:classlift}
Set $h_F:= |\Cl(O_F[1/p])|$, and suppose that
$\Phi^{\brp}_{T}(x) = 0$. Then there exists $\ti{x} \in
H^{1}_{f,\brp}(F,T_y)$ such that $j(\ti{x}) = h_F x$ and
$\Phi^{\brp}_{T_y}(\ti{x}) = 0$.
\end{proposition}

\begin{proof}
It is not hard to check that \eqref{E:ext2} yields sequences (which are
exact in the middle, and where we denote maps on resolvends by the
same symbols as the corresponding maps on cohomology groups):
\begin{align*}
&\bH(A(\bZ_p(1))) \xrightarrow{i} \bH(A(T_y)) \to \bH(A(T)), \\
&\bH(A_v(\bZ_p(1))) \xrightarrow{i_v} \bH(A_v(T_y)) \to \bH(A_v(T)), \\
&\bH(\fM(\bZ_p(1))) \xrightarrow{i} \bH(\fM(T_y)) \to \bH(\fM(T)), \\
&\bH(\fM_v(\bZ_p(1))) \xrightarrow{i_v} \bH(\fM_v(T_y)) \to \bH(\fM_v(T)).
\end{align*}
It therefore follows from \eqref{E:ali iso} that \eqref{E:ext2}
induces a sequence
\begin{equation*}
\Cl(\fM^\brp(\bZ_p(1))) \to \Cl(\fM^\brp(T_y)) \to \Cl(\fM^\brp(T))
\end{equation*}
of classgroups which is exact in the middle. We therefore deduce via
functoriality that the following diagram (whose rows are exact in the
middle) commutes:
\begin{equation} \label{E:bigclass}
\begin{CD}
\Cl(\fM^\brp(\bZ_p(1))) @>>> \Cl(\fM^\brp(T_y)) @>>> \Cl(\fM^\brp(T))
\\
@A{\Phi^{\brp}_{\bZ_p(1)}}AA   @A{\Phi^{\brp}_{T_y}}AA
@A{\Phi^\brp_T}AA \\
H^{1}_{f,\brp}(F,\bZ_p(1)) @>{i}>> H^{1}_{f,\brp}(F,T_y) @>{j}>>
H^{1}_{f,\brp}(F,T).
\end{CD}
\end{equation}

For each integer $n>0$, the Wedderburn decomposition of
$\fM^\brp(\mu_{p^n})$ (cf. \eqref{E:wedd}) is given by
$$
\fM^\brp(\mu_{p^n}) \simeq \bigoplus_{i=0}^{p^n-1} O_F[1/p].
$$
Hence $h_F \cdot \Cl(\fM^\brp(\mu_{p^n})) = 0$ for each $n$, and so it
follows that $h_F \cdot \Cl(\fM^\brp(\bZ_p(1))) = 0$ also. The result
now follows from the commutativity of \eqref{E:bigclass}.
\end{proof}

\begin{proposition} \label{P:fingreen}
Assume that the pairing $\langle \langle\,,\,\rangle \rangle$ is
non-degenerate, and suppose that $\Phi^{\brp}_{T}(x)$ is of finite
order. Then $x \in \fG_F(T)$.
\end{proposition}

\begin{proof} 
We first observe from the definition of $\fG_F(T)$ that if $M_1 x \in
\fG_F(T)$ for any integer $M_1>0$, then $x \in \fG_F(T)$ also. Hence, we
see from Proposition \ref{P:classlift} that, without loss of
generality, we may assume that there exists a lift $\ti{x} \in
H^{1}_{f,\brp}(F,T_y)$ of $x \in H^{1}_{f,\brp}(F,T)$ such that
$\Phi^{\brp}_{T_y}(\ti{x}) = 0$. We shall make this assumption from
now on.

It follows from the standard identification of $H^1(F,T^*)$ with
$\Ext^{1}_{\Omega_F}(\bZ_p,T^*)$ that we may write
\begin{equation} \label{E:glosplit}
T_y = \bZ_p(1) \times T
\end{equation}
with $\Omega_F$-action given by
$$
(\zeta, t)^\sigma = \left( \zeta^\sigma \cdot \{ f(\sigma^{-1})(t)
\}^{\sigma}, t^\sigma \right)
$$
for any fixed choice of $\Omega_F$-cocycle $f$ representing $y \in
H^1(F,T^*)$. This implies that there is an isomorphism of
$F^c$-algebras (but not of $\Omega_F$-modules)
\begin{equation} \label{E:Asplit}
A_{F^c}(T_y) \xrightarrow{\sim} A_{F^c}(\bZ_p(1)) \otimes_{F^c}
A_{F^c}(T);\qquad \alpha \mapsto k_1(\alpha) \otimes k_2(\alpha),
\end{equation}
where $k_1$ and $k_2$ are induced by the projection maps $\bZ_p(1)
\times T \to \bZ_p(1)$ and $\bZ_p(1) \times T \to T$ respectively.

We now make a choice of resolvends associated to $\ti{x}$ and
$\lambda_v$ for each $v$. Since $\ti{x} \in \Ker(\Phi^{\brp}_{T_y})$,
we may choose a resolvend $\alpha \in \bH(\fM^{\brp}(T_y))$ associated
to $\ti{x}$. For each place $v \nmid p$, we choose an arbitrary
resolvend $\nu_v \in \bH(\fM^\brp_v(T_y))$ associated to
$\lambda_v$. For each place $v \mid p$, we set
$$
\nu_v = (1 \otimes k_2(\alpha))_v
$$
(cf. \eqref{E:Asplit}); this is a resolvend associated to $\lambda_v$
since $y$ is locally trivial at $v$. Then, for each place $v$ of $F$,
we have
$$
\tau_v:= \alpha_v \nu^{-1}_{v} \in i_v(\bH(\fM^\brp_v(\bZ_p(1)))),
$$
and $\tau_v$ is a resolvend associated to $\ti{x}_v - \lambda_v \in
i_v(H^1(F_v, \bZ_p(1)))$. 

If $v \nmid p$, then $\bZ_p(1)$ is unramified at $v$, and so Remark
\ref{R:unram} implies that we have
$$
H(\fM^\brp_v(\bZ_p(1))) = H(\fM_v(\bZ_p(1))) \simeq
H^1_f(F_v,\bZ_p(1)) \simeq \check{O}^{\times}_{F_v}.
$$
Hence it follows that $(\log_p(q_v)) \ord_v(\ti{x}_v - \lambda_v) = 0$.

If $v \mid p$, then $\tau_v = (k_1(\alpha) \otimes 1)_v$, and so
$\oplus_{v \mid p} \tau_v$ is the image of $k_1(\alpha) \in
\fM^{\brp}(\bZ_p(1))^{\times}_{O_{F^c}}$ under the localisation map
$$
\delta: \fM^{\brp}(\bZ_p(1))^{\times}_{O_{F^c}} \to \bigoplus_{v \mid p}
\fM_{v}^{\brp}(\bZ_p(1))^{\times}_{O_{F^c_v}}.
$$

For each integer $n \geq 1$, write 
$$
\delta_n: \fM^\brp(\mu_{p^n})^{\times}_{O_{F^c}} \to \bigoplus_{v \mid p}
\fM^\brp_v (\mu_{p^n})^{\times}_{O_{F^c_v}}
$$
for the localisation map, and let $\tau_{v,n}$ denote the image of
$\tau_v$ in $\fM^\brp_v (\mu_{p^n})^{\times}_{O_{F^c_v}}$ under the
natural map $\fM_{v}^{\brp}(\bZ_p(1))^{\times}_{O_{F^c_v}} \to
\fM^\brp_v (\mu_{p^n})^{\times}_{O_{F^c_v}}$.  Let $(\ti{x}_v-
\lambda_v)_n$ denote the image of $(\ti{x}_v- \lambda_v)$ in
$H^1(F_v,\mu_{p^n})$ under the natural map $H^1(F_v,\bZ_p(1)) \to
H^1(F_v, \mu_{p^n})$.

Since $\bigoplus_{v \mid p} \tau_v$ lies in the image of $\delta$, we
have that $\bigoplus_{v \mid p} \tau_{v,n}$ lies in the image of
$\delta_n$. It therefore follows that 
$$
\oplus_{v \mid p}(\ti{x}_v- \lambda_v)_n \in \bigoplus_{v \mid p} H^1(F_v,
\mu_{p^n}) \simeq \bigoplus_{v \mid p} F_{v}^{\times}/ F_{v}^{\times p^n}
$$ 
lies in the image of $O_F[1/p]^{\times}F^{\times p^n}/F^{\times p^n}$
under the localisation map
$$
F^{\times}/F^{\times p^n} \to \bigoplus_{v \mid p}
F_{v}^{\times}/F_{v}^{\times p^n}.
$$
This in turn implies that 
$$
\bigoplus_{v \mid p} N_{F_{v}/\bQ_{p}} (\ti{x}_v- \lambda_v)_n \in
H^1(\bQ_p, \mu_{p^n}) \simeq \bQ_{p}^{\times}/\bQ_{p}^{\times p^n}
$$
lies in the image of $\bZ[1/p]^{\times} \bQ^{\times p^n}/\bQ^{\times
p^n}$ under the localisation map
$$
H^1(\bQ,\mu_{p^n}) \simeq \bQ^{\times}/\bQ^{\times p^n} \to H^1(\bQ_p,
\mu_{p^n}) \simeq \bQ_{p}^{\times}/\bQ_{p}^{\times p^n},
$$
for all $n \geq 1$. As the function $\log_p(z)$ vanishes on the image
of $\bZ[1/p]^\times$ in $\bQ_{p}^{\times}$, we conclude that
$$
\sum_{v \mid p} \log_p \left( N_{F_{v}/\bQ_{p}}(\ti{x}_v - \lambda_v)
\right) = 0.
$$
We therefore deduce that $B_F(x,y) = 0$. Hence $x \in \fG_F(T)$, since
by hypothesis $\langle \langle\,,\,\rangle \rangle$ is non-degenerate.
\end{proof} 

Theorem \ref{T:intro 1} now follows immediately from Propositions
\ref{P:normker} and \ref{P:fingreen}.

\end{document}